\documentclass[12pt,a4paper]{article}

\usepackage[utf8]{inputenc}
\usepackage[english]{babel}
\usepackage{amsmath,amssymb,amsfonts,amsthm,graphicx}
\usepackage{hyperref}

\begin{document}

\newcommand{\C}{\mathbb {C}}
\newcommand{\N}{\mathbb{N}}
\newcommand{\R}{\mathbb{R}}
\newcommand{\tf}{\mathcal{F}}

\renewcommand{\theenumi}{\roman{enumi}}
\renewcommand{\labelenumi}{\theenumi)}

\swapnumbers
\newtheorem{rem5}{Remark}[section]
\newtheorem{def1}[rem5]{Definition}
\newtheorem{thm1}[rem5]{Theorem}
\newtheorem{lem1}[rem5]{Lemma}
\newtheorem{lem2}[rem5]{Lemma}
\newtheorem{rem1}[rem5]{Remark}
\newtheorem{thm4}[rem5]{Theorem}
\newtheorem{rem3}[rem5]{Remark}
\newtheorem{cond1}{Condition}[section]
\newtheorem{rem4}[cond1]{Remarks}
\newtheorem{thm2}[cond1]{Theorem}
\newtheorem{rem7}[cond1]{Remark}
\newtheorem{REM0}[cond1]{Remark}
\newtheorem{thm3}[cond1]{Theorem}
\newtheorem{thm5}[cond1]{Theorem}
\newtheorem{rem6}[cond1]{Remark}
\newtheorem{cor2}[cond1]{Corollary}
\newtheorem{prop2}{Proposition}[section]
\newtheorem{prop3}[prop2]{Proposition}
\newtheorem{prop4}[prop2]{Proposition}
\newtheorem{prop5}[prop2]{Proposition}
\newtheorem{prop6}[prop2]{Proposition}
\newtheorem{cor1}[prop2]{Corollary}
\newtheorem{rem2}[prop2]{Remarks}

\title{Explicit error estimates for the stationary phase method I:\\ The influence of amplitude singularities}

\author{Felix Ali Mehmeti \footnote{Université de Valenciennes et du Hainaut-Cambrésis, LAMAV, FR CNRS 2956, Le Mont Houy, 59313 Valenciennes Cedex 9, France. Email: felix.ali-mehmeti@univ-valenciennes.fr} , Florent Dewez \footnote{Université Lille 1, Laboratoire Paul Painlevé, CNRS U.M.R 8524, 59655 Villeneuve d'Ascq Cedex, France. Email: florent.dewez@math.univ-lille1.fr}}

\date{}

\maketitle

\begin{abstract}
	We consider a version of the stationary phase method in one dimension of A. Erdélyi, allowing the phase to have stationary points of non-integer order and the amplitude to have integrable singularities. We provide a complete proof and we improve the remainder estimates in the case of regular amplitude. Then we are interested in the time-asymptotic behaviour of the solution of the free Schrödinger equation on the line, where the Fourier transform of the initial data is compactly supported and has a singularity. Applying the above mentioned method, we obtain asymptotic expansions with respect to time in certain space-time cones, where the coefficients of the remainders are uniformly bounded. These results show the influence of the singularity on the decay.
\end{abstract}

\vspace{0.3cm}



\setcounter{section}{-1}
\section{Introduction}

\hspace{2.5ex} The asymptotic behaviour of oscillatory integrals with respect to a large parameter can often be described using the stationary phase method. A theorem of A. Erdélyi \cite[section 2.8]{erdelyi} permits to treat oscillatory integrals with singular amplitudes and furnishes asymptotic expansions with explicit remainder estimates. The approach is specific for one integration variable and the results are interesting for applications. Unfortunately the proof is only sketched in the original paper \cite{erdelyi}. In the present paper, we provide a complete proof and improve the remainder estimates in the case without amplitude singularities. Then we apply the results to the solution of the free Schrödinger equation on the line for initial conditions with singular Fourier transform. We calculate asymptotic expansions to one term with respect to time as large parameter together with remainder estimates, inspired by the method of \cite{fam}. We obtain leading terms exhibiting the optimal decay rate of the solution in certain space-time cones. This result shows in which way the singularity affects the dispersion of the solution. In particular the singularity diminishes the time decay rate in certain space-time cones below the rate of quantum mechanic dispersion, when leaving the $L^2$-setting.

Consider the free Schrödinger equation
\begin{equation*}
	(S) \qquad \left\{ \begin{array}{rl}
			& \big[ i \partial_t + \partial_x^2 \big] u(t,x) = 0 \\
			& \vspace{-0.3cm} \\
			& u(0,x) = u_0(x)
	\end{array} \right. \; ,
\end{equation*}
for $t > 0$ and $x \in \R$. If we suppose $u_0 \in L^{1}(\R)$ then
\begin{equation*}
	\left\| u(t,.) \right\|_{L^{\infty}(\R)} \leqslant \frac{\| u_0 \|_{L^1(\R)}}{2\sqrt{\pi}} \, t^{-\frac{1}{2}} \; ,
\end{equation*}
see for example \cite[p.60]{reed-simon}. If it is assumed that $u_0 \in L^2(\R)$ then we have by Strichartz' estimate (\cite{strichartz}, see also \cite{banica})
\begin{equation*}
	\left\| u(t,.) \right\|_{L^{\infty}(\R)} \leqslant C \, \| u_0 \|_{L^2(\R)} \, t^{-\frac{1}{4}} \; ,
\end{equation*}
for some suitable constant $C>0$.

Consider for example initial conditions $u_0$ satisfying
\begin{equation} \label{intro}
	\forall \, p \in \R \qquad	\tf u_0 (p) = p^{\mu-1} (1-p) \, \chi_{[0,1]}(p) \; ,
\end{equation}
with $\mu \in (0,1)$; here $\tf u_0$ refers to the Fourier transform of $u_0$ and $\chi_{[0,1]}$ is the characteristic function of $[0,1]$. Under this assumption, $u_0$ is a smooth function which never belongs to $L^1(\R)$ and belongs to $L^2(\R)$ if and only if $\mu \in \big(\frac{1}{2}, 1\big)$. The question of the $L^{\infty}$-time decay rate for the above problem when $\mu \in \big(0,\frac{1}{2}\big)$ seems to be open. The theoretical results of this paper allow to make progress towards an answer in particular to this question.

In section 1, we first recall Erdélyi's result concerning asymptotic expansions with remainder estimates of oscillatory integrals of the type
\begin{equation*}
	\int_{p_1}^{p_2} U(p) \, e^{i \omega \psi(p)} \, dp \; ,
\end{equation*}
where the amplitude $U$ has singularities and the phase function $\psi$ has stationary points of real order at $p_1$ and $p_2$. We give the proof following the lines of the original demonstration: we start by splitting the integral using a cut-off function which separates the endpoints of the interval. Then we use explicit substitutions to simplify the phases. Afterwards  integrations by parts create the expansion of the integral and provide the remainder terms, that we estimate to conclude. The two last steps are carried out using complex analysis in one variable leading to explicit estimates of the error. Especially the application of the Cauchy theorem allows to shift the integration path of the integrals created by the integrations by parts into a region of controllable oscillations of the integrands.

Then we treat the case of the absence of amplitude singularities, which will be essential for certain applications. The previous estimate furnishes here only the same decay rate for the
highest term of the expansion and the remainder. The remedy proposed by Erdélyi \cite[p.55]{erdelyi} leads to complicated formulas when written down and does not seem possible in the case of stationary points of non integer order. To refine this analysis, we use the integral representation formula for the remainder obtained in section 1. Introducing a new parameter, we obtain an estimate of the integrand permitting a balance between its singular behaviour with respect to the integration variable and its decay with respect to $\omega$.

In section 2 we consider the Fourier solution formula of the free Schrödinger equation on the line with initial conditions in a compact frequency band with a singular frequency at one of the endpoints, as for example $p=0$ in \eqref{intro}. We establish an explicit remainder estimate for an asymptotic expansion of the solution formula with respect to time, using section 1 and a method based on \cite[section 3]{fam}. We seek conditions to get a uniformly bounded remainder which leads to spatial restrictions: away from the critical directions given by the endpoints of the frequency band, we obtain uniform remainder estimates in cones in space-time and we deduce from the leading term the optimal decay rate in these regions. Note that in \cite{fam}, the authors had to expand with respect to the parameter $\omega := \sqrt{t^2 + x^2}$ instead of $t$. They applied the stationary phase method as formulated in \cite{hormander}, which requires that a family of phase functions must be bounded in $\mathcal{C}^4$ in order to obtain a uniform remainder estimate. This motivated the introduction of the parameter $\omega$, which is not necessary in our context due to our explicit control of the rest in terms of the parameters $t$ and $x$. \\
In the $L^2$-case, the time decay rate $t^{-\frac{1}{2}}$ inside the cones and the fact that the spatial cross-section of any space-time cone is proportional to $t$ are used to prove that the probability amplitude behaves time-asymptotically as a laminar flow.

In the last section, we state and prove several intermediate results which have been implicitly used in \cite{erdelyi} without proof. They play a key role in Erdélyi' stationary phase method.

The above mentioned result applied to example \eqref{intro} furnishes
\begin{equation} \label{exple1}
	\left| u(t,x) - H(t,x,u_0) \, t^{-\frac{1}{2}} - K_{\mu}(t,x,u_0) \, t^{-\mu} \right| \leqslant c(u_0,\varepsilon_1,\varepsilon_2) \, \big(t^{-1} + t^{-\delta} \big) \; ,
\end{equation}
where
\begin{align}
	& \bullet \quad H(t,x,u_0) := \frac{1}{2 \sqrt{\pi}} \, e^{-i \frac{\pi}{4}} \,  e^{i \frac{x^2}{4t}} \, \left(1-\frac{x}{2t}\right) \left(\frac{x}{2t}\right)^{\mu-1} \; , \nonumber \\
	& \bullet \quad K_{\mu}(t,x,u_0) := \frac{\Gamma(\mu)}{2^{\mu+1} \pi} \, e^{i \frac{\pi \mu}{2}} \, \left(\frac{x}{2t}\right)^{-\mu} \; , \nonumber
\end{align}
for all $(t,x)$ lying in $\mathfrak{C}_{\varepsilon_1,\varepsilon_2}(0,1) := \big\{ (t,x) \, \big| \, t>0 \, , \, \varepsilon_1 \leqslant \frac{x}{2t} \leqslant 1-\varepsilon_2 \big\}$, with fixed $\varepsilon_1, \varepsilon_2 > 0$ such that $\varepsilon_1 < 1 - \varepsilon_2$ and $\delta \in \big( \max\{\mu,\frac{1}{2}\},1\big)$. In the space-time cone $\mathfrak{C}_{\varepsilon_1,\varepsilon_2}(0,1)$, the coefficients $H(t,x,u_0)$ and $K_{\mu}(t,x,u_0)$ are uniformly bounded. The expansion highlights the superposition of the quantum mechanic dispersion and the effects caused by the singularity. For $\mu < \frac{1}{2}$, the second phenomenon is dominant for large $t$. Nevertheless, the number $c(u_0, \varepsilon_1,\varepsilon_2)$ tends to infinity as $\varepsilon_1$ or $\varepsilon_2$ tends to $0$. This is connected to the fact that the above limit requires the study of an oscillatory integral changing its nature when a certain parameter attains a critical value.

We also obtain asymptotic expansions in space-time cones outside $\mathfrak{C}_{\varepsilon_1,\varepsilon_2}(0,1)$ given by $\mathfrak{C}_{1,\varepsilon}^c(0,1):= \big\{ (t,x) \, \big| \, t>0 \, , \, -\frac{1}{\varepsilon} \leqslant \frac{x}{2t} \leqslant -\varepsilon \big\}$ and $\mathfrak{C}_{2,\varepsilon}^c(0,1) := \big\{ (t,x) \, \big| \, t>0 \, , \, 1 + \varepsilon \leqslant \frac{x}{2t} \leqslant \frac{1}{\varepsilon} \big\}$ where $\varepsilon>0$ is sufficiently small; here we have
\begin{equation*}
	\Big| u(t,x) - K_{j,\mu}^c (t,x,u_0) \, t^{-\mu} \Big| \leqslant c_j^c(u_0,\varepsilon) \, t^{-1} \; ,
\end{equation*}
where for $j=1,2$
\begin{equation*}
	\begin{aligned}
		& \bullet \quad K_{j,\mu}^c(t,x,u_0) := \frac{\Gamma(\mu)}{2^{\mu+1} \pi} \, e^{(-1)^j i \frac{\pi \mu}{2}} \, \bigg( (-1)^j \Big(\frac{x}{2t}\Big)\bigg)^{-\mu} \; , \\
	\end{aligned}
\end{equation*}
for all $(t,x) \in \mathfrak{C}_{j,\varepsilon}^c(0,1)$. In this situation, the leading term is always determined by the singular frequency and the coefficients $K_{j,\mu}^c(t,x,u_0)$ are uniformly bounded in these cones. As above, a similar phenomenon produces a blow-up of $c_j^c(u_0,\varepsilon)$ when $\varepsilon$ tends to $0$.

The previous results show that in the $L^2$-case, namely $\mu \in \big(\frac{1}{2},1\big)$, wave packets in frequency bands move essentially in space-time cones. We shall prove that their $L^2$-norm inside these cones converges to a value for $t$ tending to infinity.

Indeed, we derive from \eqref{exple1} an estimate for the $L^2$-norm of the solution on the interval $I_t := [2 \varepsilon_1 \, t ,2(1-\varepsilon_2) \, t]$, which corresponds to the spatial cross-section of the cone $\mathfrak{C}_{\varepsilon_1,\varepsilon_2}(0,1)$ at time $t$,
\begin{equation*}
	\left| \left\| u(t,.) \right\|_{L^2(I_t)} - \frac{1}{\sqrt{2 \pi}} \, \left\| \tf u_0 \right\|_{L^2(\varepsilon_1, 1 - \varepsilon_2)} \right| \leqslant \tilde{c}(u_0, \varepsilon_1,\varepsilon_2) \, t^{\frac{1}{2}-\mu} \; ,
\end{equation*}
for all $t \geqslant 1$. We observe that a large part of the $L^2$-norm is concentrated in the cone in this situation, according to Plancherel's Theorem. The constant $\tilde{c}(u_0, \varepsilon_1,\varepsilon_2)$ blows up when $\varepsilon_1$ or $\varepsilon_2$ tends to $0$ since it inherits this behaviour from the above constant $c(u_0, \varepsilon_1,\varepsilon_2)$.

These observations remain true for more general initial conditions but it seems to be impossible to deduce from these results a uniform $L^{\infty}$-time decay rate because the coefficient of the remainder term of the expansion blows up when the boundaries of the cones tend to the critical directions.

Our approach in part II of this paper will be to replace the smooth cut-off function separating the stationary point and the amplitude singularity by a characteristic function. This introduces new technical difficulties but removes the artificial contribution to the blow-up of the remainder caused by the cut-off function, and permits to obtain the explicit dependence of the remainder term on the distance between the two endpoints of the interval.

The results of part I and II applied to the Schrödinger equation yield a partial answer to the above mentioned problem. As we have explained, we obtain in part I an optimal uniform estimate,
\begin{equation*}
	\left\| u(t,.) \right\|_{L^{\infty}} \leqslant C_1(u_0) \, t^{-\min\left\{\mu, \frac{1}{2} \right\}} \; ,
\end{equation*}
outside arbitrarily narrow cones centred in the critical space-time directions given by the endpoints of the frequency band. Note that it is possible to apply our version of the stationary phase method on the direction coming from the singular frequency. In this case, we obtain
\begin{equation*}
	\left| u(t,x) - L_{\mu}(t,u_0) \, t^{-\frac{\mu}{2}} \right| \leqslant c(u_0) \left( t^{-1} + t^{-\frac{1}{2}} \right) \; ,
\end{equation*}
where the coefficient $L_{\mu}(t,u_0)$ is explicitly given (see \eqref{Lmu}).
In part II, we shall replace the un-controllable cone generated by the singularity by a smaller region delimited by a family of curves whose narrowness is parametrized by some positive $\varepsilon$. Outside this region and between the two critical directions, we obtain an optimal uniform estimate
\begin{equation*}
	\left\| u(t,.) \right\|_{L^{\infty}} \leqslant C_2(u_0, \varepsilon) \, t^{-\frac{\mu}{2} - \varepsilon} \; ,
\end{equation*}
but still with the blow-up of $C_2(u_0, \varepsilon)$ when $\varepsilon$ tends to $0$. These results lead us to the conjecture that a global optimal estimate
\begin{equation*}
	\left\| u(t,.) \right\|_{L^{\infty}(\R)} \leqslant C_3(u_0) \, t^{-\frac{\mu}{2}}
\end{equation*}
holds.


Finally, we comment on some related results. The time-decay rate of the free Schrö\-din\-ger equation is considered in \cite{cazenave1998} and \cite{cazenave2010}. In \cite{cazenave1998}, singular initial conditions are constructed to derive the exact $L^p$-time decay rates of the solution, which are slower than the classical results for regular initial conditions. In \cite{cazenave2010}, the authors construct initial conditions in Sobolev spaces (based on the Gaussian function), and they show that the related solutions has no definite $L^p$-time decay rates, nor coefficients, even though upper estimates for the decay rates are established.\\
The papers \cite{cazenave1998} and \cite{cazenave2010} use special formulas for functions and their Fourier transforms, which are themselves based on complex analysis. In our results, we furnish slower decay rates by considering initial conditions with singular Fourier transforms. Here, complex analysis is directly applied to the solution formula of the equation, which permits to obtain results for a whole class of functions. The method seems to be more flexible.

In \cite{fam0}, the authors study the time-asymptotic behaviour of the solution of the Schrö\-din\-ger equation on star-shaped networks with a localized potential. They establish a perturbation inequality which shows that the evolution of high frequency signals is close to the evolution of the free equation. This points out the usefulness of detailed informations on the motion of wave packets in frequency bands.

The article \cite{fam2} considers the Klein-Gordon equation on $\R$ with potential step. The authors introduce the idea of considering frequency bands to obtain qualitative informations on the solution.

In the setting of the previous paper, the aim of \cite{fam3} is to describe the influence of the height of a potential step on the time-asymptotic energy flow of wave packets. To do so, the authors use asymptotic expansions of the solution in certain space-time cones (as in \cite{fam}) and need an extra assumption on the initial data to refine an estimate from below.

In \cite{banica}, the author considers one-dimensional Schrödinger equations with singular coefficients. Dispersion inequalities and Strichartz-type estimates are furnished.

In \cite[chapter 4]{evans}, a stationary phase method is provided. The author assumes that the amplitude belongs to $\mathcal{C}_c^{\infty}(\R^d)$ (for $d \geqslant 1$) and supposes that the stationary points of the phase are non-degenerate. Firstly the author employs Morse's lemma to simplify the phase function. Then using Fubini's theorem, he obtains a product of functions, where each of them is a Fourier transform of a tempered distribution. Finally computing these Fourier transforms and estimating them leads to the result.
Nevertheless, the use of Morse's lemma implies a loss of precision regarding the estimate of the remainder. Indeed, Morse's lemma is based on the implicit function theorem and so the substitution is not explicit.

The paper \cite{hormander2} treats the propagation of singularities in space-time for operators with real characteristics. Note that in the case of the hyperbolic problem \cite{fam}, the space-time cones tend to the characteristic set if the frequency band tends to infinity.

In \cite[chapter 7]{hormander}, the author provides a stationary phase method which is different from \cite{evans}. It is assumed that the amplitude $U$ and the phase have a certain regularity on $\R^d$ (for $d \geqslant 1$) and that $U$ has a compact support. Asymptotic expansions of the oscillatory integral are given by using Taylor's formula of the phase, where the stationary point is supposed non-degenerate. Morse's lemma is not needed in this situation. However, stronger hypothesis concerning the phase are required in order to bound uniformly the remainder by a constant, which is not explicit.

The article \cite{liess} deals with the time decay rates for the system of crystal optics. The stationary phase method \cite{hormander} is employed to obtain the decay rates. Observe that a change of parameters is carried out to obtain a bounded phase in $\mathcal{C}^4$, which permits to apply this stationary phase method, like in \cite{fam}.

In \cite{strichartz}, the author furnishes the first Strichartz-type estimates for the Schrödinger equation and the Klein-Gordon equation. Using complex analysis, the author provides estimates of the $L^2(S)$-norm of Fourier transforms of functions belonging to $L^q(\R^d)$, for some $q \geqslant 1$, where $S$ is a quadratic surface. These considerations lead to the above mentioned estimates.\\

\noindent \textbf{Acknowledgements:}\\
The authors thank S. De Bièvre, R. Haller-Dintelmann and V. Régnier for valuable discussions.

\section{Stationary points of real order and singular amplitudes: explicit error estimates in one variable}

\hspace{2.5ex} Before formulating the results of this section, let us introduce the two assumptions related to the amplitude and to the phase.

Let $p_1,p_2$ be two real numbers such that $-\infty < p_1 < p_2 < +\infty$.\\ \\
\textbf{Assumption (P$_{\rho_1,\rho_2,N}$).} For $\rho_1, \rho_2 \geqslant 1$, let $\psi \in \mathcal{C}^1\big([p_1,p_2], \R\big)$ be a function satisfying
	\begin{equation*}
		\forall \, p \in [p_1,p_2] \qquad \psi'(p) = (p-p_1)^{\rho_1 -1} (p_2-p)^{\rho_2 -1} \tilde{\psi}(p) \; ,
	\end{equation*}
	where $\tilde{\psi} \in \mathcal{C}^N\big([p_1,p_2], \R\big)$ is assumed positive, with $N \in \N \backslash \{0\}$. Points $p_j$ $(j=1,2)$ are called \emph{stationary points} of $\psi$ of order $\rho_j -1$.\\
	
\noindent \textbf{Assumption (A$_{\mu_1,\mu_2,N}$).} For $0 < \mu_1 , \mu_2 \leqslant 1$, let $U : (p_1, p_2) \longrightarrow \C$ be a function defined by
	\begin{equation*}
		\forall \, p \in (p_1,p_2) \qquad U(p )= (p - p_1)^{\mu_1 -1} (p_2 - p)^{\mu_2 -1} \tilde{u}(p) \; ,
	\end{equation*}
	where $\tilde{u} \in \mathcal{C}^N\big([p_1,p_2], \C\big)$, with $N \in \N \backslash \{0\}$, and $\tilde{u}(p_j) \neq 0$ if $\mu_j \neq 1$ $(j=1,2)$. Points $p_j$ are called \emph{singularities} of $U$.

\begin{rem5}
	\em The hypothesis $\tilde{u}(p_j) \neq 0$ if $\mu_j \neq 1$ prevents the function $\tilde{u}$ from affecting the behaviour of the singularity $p_j$.
\end{rem5}

\subsection*{Non-vanishing singularities: Erdelyi’s theorem}

\hspace{2.5ex} The aim of this subsection is to state Erdélyi's result \cite[section 2.8]{erdelyi} and give a complete proof.

Let us define some objects that will be used throughout the rest of this paper.

\begin{def1} \label{DEF1}
	\begin{enumerate}
	\item For $j=1,2$, let $\varphi_j : I_j \longrightarrow \R$ be the functions defined by
	\begin{equation*}
		\varphi_1(p) := \big( \psi(p) - \psi(p_1) \big)^{1/\rho_1} \qquad \text{and} \qquad \varphi_2(p) := \big( \psi(p_2) - \psi(p) \big)^{1/\rho_2} \; ,
	\end{equation*}
	with $I_1 := [p_1, p_2 - \eta]$, $I_2 := [p_1 + \eta, p_2]$ and $s_1 := \varphi_1(p_2 - \eta)$, $s_2 := \varphi_2(p_1 + \eta)$.
	\item For $j=1,2$, define $k_j : (0,s_j] \longrightarrow \C$ by
	\begin{equation*}
			k_j(s) := U\big( \varphi_j^{-1}(s) \big) \, s^{1- \mu_j} \, \big( \varphi_j^{-1} \big)'(s) \; ,
	\end{equation*}
	which can be extended to the interval $[0,s_j]$ (see Proposition \ref{PROP4}).
	\item Let $\nu : [p_1,p_2] \longrightarrow \R$ be a smooth function such that
	\begin{equation*}
		\left\{ \begin{array}{rl}
					& \nu = 1 \qquad \text{on} \quad [p_1, p_1 + \eta] \; , \\
					& \nu = 0 \qquad \text{on} \quad [p_2 - \eta, p_2] \; ,\\
					& 0 \leqslant \nu \leqslant 1 \; ,
				\end{array} \right.
	\end{equation*}
	with fixed $\eta \in \big(0, \frac{p_2 - p_1}{2}\big)$. For $j=1,2$, define $\nu_j : [0,s_j] \longrightarrow \R$ by
	\begin{equation*}
		\nu_1(s) := \nu \circ \varphi_1^{-1}(s) \qquad \text{and} \qquad \nu_2(s) := (1 - \nu) \circ \varphi_2^{-1}(s) \; .
	\end{equation*}
	\item For $s > 0$, the complex curve $\Lambda^{(j)}(s)$ is defined by
	\begin{equation*}
		\Lambda^{(j)}(s) := \left\{ s + t e^{(-1)^{j+1} i \frac{\pi}{2\rho_j}} \, \Big| \, t \geqslant 0 \right\} \; .
	\end{equation*}
	\end{enumerate}
\end{def1}

\vspace{0.2cm}
	
\begin{thm1} \label{THM1}
	Let $N \in \N \backslash \{0\}$, let $\rho_1, \rho_2 \geqslant 1$ and $0 < \mu_1, \mu_2 < 1$. Suppose that the functions $\psi : [p_1,p_2] \longrightarrow \R$ and $U : (p_1, p_2) \longrightarrow \C$ satisfy Assumption \emph{(P$_{\rho_1,\rho_2,N}$)} and Assumption \emph{(A$_{\mu_1,\mu_2,N}$)}, respectively. Then there exist functions $A_N^{(j)}$, $R_N^{(j)} : (0,+\infty) \longrightarrow \C$ for $j=1,2$, such that:
	\begin{equation*}
		\left\{ \begin{array}{rl}
			& \displaystyle \int_{p_1}^{p_2} U(p) e^{i \omega \psi(p)} \, dp = \sum_{j=1,2} \left( A_N^{(j)}(\omega) + R_N^{(j)}(\omega) \right) \; , \\
			& \displaystyle \left| R_N^{(j)}(\omega) \right| \leqslant \frac{1}{(N-1)!} \, \frac{1}{\rho_j} \, \Gamma\bigg(\frac{N}{\rho_j}\bigg) \int_0^{s_j} s^{\mu_j-1} \left| \frac{d^N}{ds^N}(\nu_j k_j)(s) \right| \, ds \; \omega^{-\frac{N}{\rho_j}} \; ,
		\end{array} \right.
	\end{equation*}
	for all $\omega > 0$. For $j =1,2$ and $\omega > 0$, we have defined
	\begin{equation*}
		\begin{aligned}
			\bullet \quad & A_N^{(j)}(\omega) := e^{i \omega \psi(p_j)} \, \sum_{n=0}^{N-1} \Theta_{n+1}^{(j)}(\rho_j,\mu_j) \, \frac{d^n}{ds^n}(k_j)(0) \, \omega^{-\frac{n+\mu_j}{\rho_j}} \; , \\
			\bullet \quad & R_N^{(j)}(\omega) := (-1)^{N+1+j} \, e^{i\omega \psi(p_j)} \int_0^{s_j} \phi_{N}^{(j)}(s,\omega, \rho_j, \mu_j) \, \frac{d^N}{ds^N}(\nu_j k_j)(s) \, ds \; ,
		\end{aligned}
	\end{equation*}
	where, for $n=0,\ldots,N-1$:
	\begin{equation*}
		\begin{aligned}
			\bullet \quad & \Theta_{n+1}^{(j)}(\rho_j, \mu_j) := \frac{(-1)^{j+1}}{n! \, \rho_j} \, \Gamma \Big(\frac{n+\mu_j}{\rho_j}\Big) \, e^{(-1)^{j+1} i \frac{\pi}{2} \, \frac{n+\mu_j}{\rho_j}} \; , \\
			\bullet \quad & \phi_{n+1}^{(j)}(s,\omega,\rho_j, \mu_j) := \frac{(-1)^{n+1}}{n!} \int_{\Lambda^{(j)}(s)} (z-s)^n z^{\mu_j -1} e^{(-1)^{j+1} i\omega z^{\rho_j}} \, dz \; .
		\end{aligned}
	\end{equation*}
\end{thm1}

\begin{proof}
	For fixed $\rho_j \geqslant 1$ and $0 < \mu_j < 1$, we shall note $\phi_{n}^{(j)}(s,\omega)$ instead of $\phi_{n}^{(j)}(s,\omega,\rho_j,\mu_j)$. Now we take $\omega > 0$ and we divide the proof in five steps.\\
	
	\noindent \textit{First step: Splitting of the integral.} Using the cut-off function $\nu$, we can write the integral as follows:
	\begin{equation*}
		\int_{p_1}^{p_2} U(p) e^{i \omega \psi(p)} \, dp = \tilde{I}^{(1)}(\omega) + \tilde{I}^{(2)}(\omega) \; ,
	\end{equation*}
	where
	\begin{equation*}
		\tilde{I}^{(1)}(\omega) := \int_{p_1}^{p_2-\eta} \nu(p) U(p) \, e^{i \omega \psi(p)} \, dp \quad , \quad \tilde{I}^{(2)}(\omega) := \int_{p_1+\eta}^{p_2} \big( 1-\nu(p) \big) U(p) \, e^{i \omega \psi(p)} \, dp \; .
	\end{equation*}
	\textit{Second step: Substitution.} Proposition \ref{PROP3} shows that $\varphi_j : I_j \longrightarrow [0,s_j]$ is a $\mathcal{C}^{N+1}$-diffeomorphism. Using the substitution $s = \varphi_1(p)$, we get
	\begin{equation*}
		\begin{aligned}
			\tilde{I}^{(1)}(\omega)	& = \int_{p_1}^{p_2-\eta} \nu(p) U(p) \, e^{i \omega \psi(p)} \, dp \\
								& = e^{i \omega \psi(p_1)} \int_0^{s_1} \nu \big( \varphi_1^{-1}(s) \big) U \big( \varphi_1^{-1} \big) \, e^{i \omega s^{\rho_1}} (\varphi_1^{-1})'(s) \, ds \\
								& = e^{i \omega \psi(p_1)} \int_0^{s_1} \nu \big( \varphi_1^{-1}(s) \big) U \big( \varphi_1^{-1} \big) s^{1-\mu_1} (\varphi_1^{-1})'(s) \, s^{\mu_1-1}  e^{i \omega s^{\rho_1}} \, ds \\
								& = e^{i \omega \psi(p_1)} \int_0^{s_1} \nu_1(s) k_1(s) \, s^{\mu_1-1} e^{i \omega s^{\rho_1}} \, ds \; ,
		\end{aligned}
	\end{equation*}
	where $k_1$ and $\nu_1$ are introduced in Definition \ref{DEF1}. In a similar way, we obtain
	\begin{equation*}
		\tilde{I}^{(2)}(\omega) = - e^{i \omega \psi(p_2)} \int_0^{s_2} \nu_2(s) k_2(s) \, s^{\mu_2-1} e^{-i \omega s^{\rho_2}} \, ds \; .
	\end{equation*}

	\vspace{0.3cm}

	\noindent \textit{Third step: Integrations by parts.} Corollary \ref{COR1} provides successive primitives of the function $s \mapsto s^{\mu_j - 1} e^{(-1)^{j+1} i\omega s^{\rho_j}}$. Moreover Proposition \ref{PROP4} ensures that $k_j \in \mathcal{C}^N([0,s_j])$. Thus by $N$ integrations by parts, we obtain
	\begin{equation*}
		\begin{aligned}
			e^{-i \omega \psi(p_1)} \tilde{I}^{(1)}(\omega)	& = \int_0^{s_1} \nu_1(s) k_1(s) \, s^{\mu_1-1} e^{i \omega s^{\rho_1}} \, ds \\
														& = \Big[ \phi_1^{(1)}(s,\omega) \big( \nu_1 k_1 \big)(s)\Big]_0^{s_1} - \int_0^{s_1} \phi_1^{(1)}(s,\omega) \, \frac{d}{ds}\big( \nu_1 k_1 \big)(s) \, ds \\
														& = \dots \\
														& = \sum_{n=0}^{N-1} (-1)^n \Big[ \phi_{n+1}^{(1)}(s,\omega) \, \frac{d^n}{ds^n}\big( \nu_1 k_1 \big)(s)\Big]_0^{s_1} \\
														& \qquad \qquad + (-1)^N \int_0^{s_1} \phi_{N}^{(1)}(s,\omega) \, \frac{d^N}{ds^N}\big( \nu_1 k_1 \big)(s) \, ds \; .
		\end{aligned}
	\end{equation*}
	One can simplify the last expression using the properties of the function $\nu_1$; indeed by hypothesis, we have $\nu(p_1)=1$, $\nu(p_2-\eta)=0$ and $\displaystyle \frac{d^n}{dp^n}(\nu)(p_1)=\frac{d^n}{dp^n}(\nu)(p_2-\eta)=0$, for $n \geqslant 1$. So the definition of $\nu_1$ implies
	\begin{equation*}
		\nu_1(0) = \nu(p_1) = 1 \qquad \text{and} \qquad \nu_1(s_1) = \nu(p_2-\eta) = 0 \; ,
	\end{equation*}
	and by the product rule applied to $\nu_1 k_1$, it follows,
	\begin{equation*}
		\frac{d^n}{ds^n}\big( \nu_1 k_1 \big)(0) = \frac{d^n}{ds^n}(k_1)(0) \qquad \text{and} \qquad \frac{d^n}{ds^n}\big( \nu_1 k_1 \big)(s_1) = 0 \; .
	\end{equation*}
	Finally,
	\begin{equation*}
		\begin{aligned}
			\tilde{I}^{(1)}(\omega)	& = \sum_{n=0}^{N-1} (-1)^{n+1} \phi_{n+1}^{(1)}(0,\omega) \, \frac{d^n}{ds^n} (k_1)(0) \, e^{i \omega \psi(p_1)}\\
							& \qquad \qquad + (-1)^{N} \, e^{i \omega \psi(p_1)} \int_0^{s_1} \phi_{N}^{(1)}(s,\omega) \, \frac{d^N}{ds^N}\big( \nu_1 k_1 \big)(s) \, ds \; .
		\end{aligned}
	\end{equation*}
	In a similar way, we have
	\begin{equation*}
		\begin{aligned}
			\tilde{I}^{(2)}(\omega)	& = \sum_{n=0}^{N-1} (-1)^n \phi_{n+1}^{(2)}(0,\omega) \, \frac{d^n}{ds^n} (k_2)(0) \, e^{i \omega \psi(p_2)}\\
							& \qquad \qquad + (-1)^{N+1} \, e^{i \omega \psi(p_2)} \int_0^{s_2} \phi_{N}^{(2)}(s,\omega) \, \frac{d^N}{ds^N}\big( \nu_2 k_2 \big)(s) \, ds \; . \\
		\end{aligned}
	\end{equation*}
	
	\vspace{0.4cm}
	
	\noindent \textit{Fourth step: Calculation of the main terms.} Let us make clear the main terms of $\tilde{I}^{(j)}(\omega)$ by calculating the coefficient $\phi_{n+1}^{(j)}(0,\omega)$, which is well-defined (see Remarks \ref{REM2}). Recall the expression of $\phi_{n+1}^{(j)}(s,\omega)$ from Corollary \ref{COR1}:
	\begin{equation*}
		\phi_{n+1}^{(j)}(s,\omega) = \frac{(-1)^{n+1}}{n!} \int_{\Lambda^{(j)}(s)} (z-s)^n \, z^{\mu_j - 1} \, e^{(-1)^{j+1} i \omega z^{\rho_j}} \, dz \; ,
	\end{equation*}
	for all $s \in [0,s_j]$ and $n=0,\ldots,N-1$. Choose $j=1$, put $s=0$ and parametrize the curve $\Lambda^{(1)}(0)$ with $z = t e^{i \frac{\pi}{2 \rho_1}}$ for $t \geqslant 0$; this leads to
	\begin{equation*}
		\phi_{n+1}^{(1)}(0,\omega) = \frac{(-1)^{n+1}}{n!} \, e^{i \frac{\pi}{2} \, \frac{n+\mu_1}{\rho_1}} \int_0^{+\infty} t^{n+\mu_1 - 1} \, e^{-\omega t^{\rho_1}} \, dt \; .
	\end{equation*}
	Setting $y= \omega t^{\rho_1}$ in the previous integral gives
	\begin{equation*}
		\begin{aligned}
			\phi_{n+1}^{(1)}(0,\omega)	& = \frac{(-1)^{n+1}}{n!} \, e^{i \frac{\pi}{2} \, \frac{n+\mu_1}{\rho_1}} \, (\rho_1 \, \omega)^{-1} \int_0^{+\infty} \Big(\frac{y}{\omega} \Big)^{\frac{n+\mu_1}{\rho_1} - 1} \, e^{-y} \, dy \\
										& = \frac{(-1)^{n+1}}{n!} \, e^{i \frac{\pi}{2} \, \frac{n+\mu_1}{\rho_1}} \, \frac{1}{\rho_1} \, \Gamma\Big(\frac{n + \mu_1}{\rho_1} \Big) \, \omega^{-\frac{n+\mu_1}{\rho_1}} \; ,
		\end{aligned}
	\end{equation*}
	where $\Gamma$ is the Gamma function defined by
	\begin{equation*}
		\Gamma \, : \, z \in \big\{ z \in \C \, \big| \, \Re(z) > 0 \big\} \, \longmapsto \, \int_0^{+\infty} t^{z-1} e^{-t} \, dt \in \C \; .
	\end{equation*}
	A similar work provides
	\begin{equation*}
		\phi_{n+1}^{(2)}(0,\omega) = \frac{(-1)^{n+1}}{n!} \, e^{-i \frac{\pi}{2} \, \frac{n+\mu_2}{\rho_2}} \, \frac{1}{\rho_2} \, \Gamma\Big(\frac{n + \mu_2}{\rho_2} \Big) \, \omega^{-\frac{n+\mu_2}{\rho_2}} \; .
	\end{equation*}
	Then we obtain
	\begin{equation*}
		e^{i \omega \psi(p_j)} \sum_{n=0}^{N-1} (-1)^{n+j} \phi_{n+1}^{(j)}(0,\omega) \, \frac{d^n}{ds^n}(k_j)(0) = e^{i \omega \psi(p_j)} \sum_{n=0}^{N-1} \Theta_{n+1}^{(j)}(\rho_j,\mu_j) \, \frac{d^n}{ds^n}(k_j)(0) \, \omega^{- \frac{n+\mu_j}{\rho_j} } \; ,
	\end{equation*}
	where we set : $\displaystyle \Theta_{n+1}^{(j)}(\rho_j,\mu_j) := \frac{(-1)^{j+1}}{n! \, \rho_j} \, \Gamma \Big( \frac{n + \mu_j}{\rho_j} \Big) \, e^{(-1)^{j+1} \frac{\pi}{2} \, \frac{n + \mu_j}{\rho_j}}$.
	
	\vspace{0.4cm}
	
	\noindent \textit{Fifth step: Remainder estimates.} The last step consists in estimating the remainders $R_N^{(j)}(\omega)$. For $j=1,2$, one can see for all $s \in (0,s_j]$ and for all $t \geqslant 0$,
	\begin{equation}\label{est3}
		s \leqslant \left| s + t e^{(-1)^{j+1} i \frac{\pi}{2 \rho_j}} \right| \qquad \Longrightarrow \qquad s^{\mu_j-1} \geqslant \left| s + t e^{(-1)^{j+1} i \frac{\pi}{2 \rho_j}} \right|^{\mu_j-1} \; ,
	\end{equation}
	since $\mu_j \in (0,1)$. Use Proposition \ref{PROP2}, inequality \eqref{est3} and parametrize $\phi_{N}(s,\omega)$ with $z = s + t e^{(-1)^{j+1} i \frac{\pi}{2 \rho_j}}$ for $t \geqslant 0$ to obtain
	\begin{align}
		\Big| \, \phi_{N}^{(j)}(s,\omega) \Big|	& \leqslant \frac{1}{(N-1)!} \int_0^{+\infty} t^{N-1} \left| \, s + t e^{(-1)^{j+1} i\frac{\pi}{2 \rho_j}} \right|^{\mu_j -1} \left| \, e^{(-1)^{j+1} i \omega \big( s + t e^{(-1)^{j+1} i\frac{\pi}{2 \rho_j}} \big)^{\rho_j}} \right| \, dt \nonumber \\
													& \leqslant \frac{1}{(N-1)!} \, s^{\mu_j-1} \int_0^{+\infty} t^{N-1} e^{-\omega t^{\rho_j}} \, dt \nonumber \\
													& = \label{est4} \frac{1}{(N-1)!} \, s^{\mu_j-1} \, \frac{1}{\rho_j} \, \Gamma\left(\frac{N}{\rho_j}\right) \omega^{-\frac{N}{\rho_j}} \; ,
	\end{align}
	where the last equality is obtained by using the substitution $y = \omega \, t^{\rho_j}$. Employing the definition of $R_N^{(j)}(\omega)$ and inequality \eqref{est4} leads to
	\begin{equation}\label{rest1}
		\begin{aligned}
			\Big| \, R_{N}^{(j)}(\omega) \Big|	& \leqslant \int_0^{s_j} \left| \, \phi_{N}^{(j)}(s,\omega,\rho_j, \mu_j) \right| \, \left| \frac{d^N}{ds^N}(\nu_j k_j)(s) \right| \, ds \\
												& \leqslant \frac{1}{(N-1)!} \, \frac{1}{\rho_j} \, \Gamma\bigg(\frac{N}{\rho_j}\bigg) \int_0^{s_j} s^{\mu_j-1} \left| \frac{d^N}{ds^N}(\nu_j k_j)(s) \right| \, ds \; \omega^{-\frac{N}{\rho_j}} \; .
		\end{aligned}
	\end{equation}
	We note that the last integral is well-defined because $\displaystyle \frac{d^N}{ds^N}(\nu_j k_j) : [0,s_j] \longrightarrow \R$ is continuous and $s \longmapsto s^{\mu_j - 1}$ is locally integrable on $[0,s_j]$.\\
	We remark finally that the decay rates of $A_N^{(j)}(\omega)$ and $R_N^{(j)}(\omega)$ are $\omega^{-\frac{N-1+\mu_j}{\rho_j}}$ and $\omega^{-\frac{N}{\rho_j}}$, respectively. Thus we observe that the decay rate of the remainder with respect to $\omega$ is higher than the one of the highest term of the expansion. This ends the proof.
\end{proof}

\subsection*{Amplitudes without singularities}

\hspace{2.5ex} The preceding theorem remains true if we suppose $\mu_j = 1$, that is to say the amplitude $U$ is regular at the point $p_j$. But in this case, we observe that the decay rate of the  highest term of the expansion and the remainder are the same, namely $\omega^{-\frac{N}{\rho_j}}$. The aim of this subsection is to refine the estimate of the remainder in this specific case.

For this purpose, the two following lemmas are used. In the first one, we provide two estimates of the function $s \longmapsto \phi_{N}^{(j)}(s,\omega,\rho_j,1)$, an estimate has a right-hand side which is singular with respect to the variable $s$ whereas the other one is regular. We carry out integrations by parts to establish this result.

\begin{lem1} \label{LEM1}
	Fix $j \in \{1,2\}$, $\rho_j \geqslant 1$ and $N \in \N \backslash \{0\}$. Then there exist three constants $a_{N,\rho_j}, b_{N,\rho_j}, c_{N,\rho_j} > 0$ such that
	\begin{equation*}
		\left\{ \begin{array}{rl}
			& \displaystyle \left| \phi_{N}^{(j)}(s,\omega,\rho_j,1) \right| \leqslant a_{N,\rho_j} \, \omega^{-\frac{N}{\rho_j}} \; , \\
			& \displaystyle \left| \phi_{N}^{(j)}(s,\omega,\rho_j,1) \right| \leqslant b_{N,\rho_j} \, \omega^{-\left(1+\frac{N-1}{\rho_j}\right)} \, s^{1-\rho_j} \, + c_{N,\rho_j} \, \omega^{-\left(1+\frac{N}{\rho}\right)} \, s^{-\rho_j} \; ,
		\end{array} \right.
	\end{equation*}
	for all $s, \omega > 0$. The constants $a_{N,\rho_j}, b_{N,\rho_j}, c_{N,\rho_j} > 0$ are given in the proof.	
\end{lem1}

\begin{rem1}
	\em Note that one can extend $\phi_{N}^{(j)}(.,\omega,\rho_j,1) \, : \, (0,s_j] \longrightarrow \R$ to $(0,+\infty)$, see Remarks \ref{REM2}.
\end{rem1}

\begin{proof}
	Fix $s > 0$, $\omega > 0$ and choose for example $j=1$. Recall the expression of $\phi_{N}^{(1)}(s,\omega,\rho_1,1)$ with the parametrization of the path $\Lambda^{(1)}(s)$ given in Definition \ref{DEF1}:
	\begin{equation*}
		\phi_{N}^{(1)}(s,\omega,\rho_1,1) = \frac{(-1)^N}{(N-1)!} \int_0^{+\infty} t^{N-1} e^{i \frac{\pi (N-1)}{2 \rho_1}} \, e^{i \omega \left(s+te^{i \frac{\pi}{2 \rho_1}}\right)^{\rho_1}} \, dt \, e^{i \frac{\pi}{2 \rho_1}} \; .
	\end{equation*}
	On the one hand, estimate \eqref{est4} is still valid for $\mu_1 = 1$:
	\begin{equation*}
		\left| \phi_{N}^{(1)}(s,\omega,\rho_1,1) \right| \leqslant \frac{1}{(N-1)!} \, \frac{1}{\rho_1} \, \Gamma\left(\frac{N}{\rho_1}\right) \, \omega^{-\frac{N}{\rho_1}} \; .
	\end{equation*}
	Put $a_{N,\rho_1} := \frac{1}{(N-1)!} \, \frac{1}{\rho_1} \, \Gamma\Big(\frac{N}{\rho_1}\Big)$ then we get the first estimate of the lemma.\\
	On the other hand, we establish the second inequality by using integrations by parts. To do so, remark that for all $s>0$ the first derivative of the function $t \in (0,+\infty) \longmapsto i \omega \big(s+te^{i \frac{\pi}{2 \rho_1}}\big)^{\rho_1}$ does not vanish; therefore we can write
	\begin{equation*}
		e^{i \omega \left(s+te^{i \frac{\pi}{2 \rho_1}}\right)^{\rho_1}} = (i \omega \rho_1)^{-1} e^{-i \frac{\pi}{2 \rho_1}} \left(s+te^{i \frac{\pi}{2 \rho_1}}\right)^{1-\rho_1} \, \frac{d}{dt}\left[e^{i \omega \left(s+te^{i \frac{\pi}{2 \rho_1}}\right)^{\rho_1}}\right] \; .
	\end{equation*}
	Moreover Proposition \ref{PROP2} implies
	\begin{equation} \label{estinf}
		\forall \, s > 0 \qquad \bigg| e^{i \omega \left(s+te^{i \frac{\pi}{2 \rho_1}}\right)^{\rho_1}} \bigg| \leqslant e^{- \omega t^{\rho_1}} \longrightarrow  0 \quad , \quad t \longrightarrow +\infty \; .
	\end{equation}
	Now we distinguish the two following cases:
	\begin{itemize}
		\item \textit{Case $N=1$.} Thanks to the two previous observations, we obtain
		\begin{equation*}
			\begin{aligned}
				\phi_1^{(1)}(s,\omega,\rho_1,1)	& = - (i \omega \rho_1)^{-1} \int_0^{+\infty} \left(s+te^{i \frac{\pi}{2 \rho_1}}\right)^{1-\rho_1} \frac{d}{dt}\left[e^{i \omega \left(s+te^{i \frac{\pi}{2 \rho_1}}\right)^{\rho_1}}\right] \, dt \\
													& = (i \omega \rho_1)^{-1} s^{1-\rho_1} e^{i \omega s^{\rho_1}} \\
													& \qquad + \frac{1-\rho_1}{i \omega \rho_1} \, e^{i \frac{\pi}{2 \rho_1}} \int_0^{+\infty} \left(s+te^{i \frac{\pi}{2 \rho_1}}\right)^{-\rho_1} e^{i \omega \left(s+te^{i \frac{\pi}{2 \rho_1}}\right)^{\rho_1}} \, dt \; ,
			\end{aligned}
		\end{equation*}
		where the boundary term at infinity is zero according to \eqref{estinf}. It follows:
		\begin{align}
			\Big| \phi_1^{(1)}(s,\omega,\rho_1,1) \Big|	& \leqslant (\omega \rho_1)^{-1} s^{1-\rho_1} \nonumber \\
															& \qquad + \frac{\rho_1-1}{\omega \rho_1} \int_0^{+\infty} \left| s+te^{i \frac{\pi}{2 \rho_1}}\right|^{-\rho_1} \left|e^{i \omega \left(s+te^{i \frac{\pi}{2 \rho_1}}\right)^{\rho_1}}\right| \, ds \nonumber \\
															& \label{ineq1}\leqslant (\omega \rho_1)^{-1} s^{1-\rho_1} + \frac{\rho_1-1}{\omega \rho_1} \, s^{-\rho_1} \int_0^{+\infty} e^{-\omega t^{\rho_1}} \, dt \\
															& = \label{eq1} \frac{1}{\rho_1} \, \omega^{-1} s^{1-\rho_1} + \frac{\rho_1 -1}{\rho_1^2} \, \Gamma\left(\frac{1}{\rho_1}\right) \, \omega^{-\left(1+\frac{1}{\rho_1}\right)} s^{-\rho_1} \; ;
		\end{align}
		Proposition \ref{PROP2} permits to obtain \eqref{ineq1} by estimating the exponential and we use the substitution $y=\omega t^{\rho_1}$ to get \eqref{eq1}. Then we set $b_{1,\rho_1} := \rho_1^{-1}$ and $c_{1,\rho_1} := \frac{\rho_1 -1}{\rho_1^2} \Gamma\Big(\frac{1}{\rho_1}\Big)$.\\
		\item \textit{Case $N \geqslant 2$.} We proceed as above by using an integration by parts:
		\begin{align}
			\phi_{N}^{(1)}(s,\omega,\rho_1,1)	& = \frac{(-1)^N}{(N-1)!} \, e^{i \frac{\pi (N-1)}{2 \rho_1}} (i \omega \rho_1)^{-1} \nonumber \\
												& \qquad \times \int_0^{+\infty} t^{N-1} \left(s+te^{i \frac{\pi}{2 \rho_1}}\right)^{1-\rho_1} \frac{d}{dt} \left[ e^{i \omega \left(s+te^{i \frac{\pi}{2 \rho_1}}\right)^{\rho_1}} \right] \, dt \nonumber \\
												& = \frac{(-1)^{N+1}}{(N-1)!} \, e^{i \frac{\pi (N-1)}{2 \rho_1}} (i \omega \rho_1)^{-1} \nonumber \\
												& \label{ipp} \qquad \times \int_0^{+\infty} \frac{d}{dt} \left[ t^{N-1} \left(s+te^{i \frac{\pi}{2 \rho_1}}\right)^{1-\rho_1} \right] e^{i \omega \left(s+te^{i \frac{\pi}{2 \rho_1}}\right)^{\rho_1}} \, dt \\
												& = \frac{(-1)^{N+1}}{(N-1)!} \, e^{i \frac{\pi (N-1)}{2 \rho_1}} (i \omega \rho_1)^{-1} \nonumber \\
												& \qquad \times \bigg( (N-1) \int_0^{+\infty} t^{N-2} \left(s+te^{i \frac{\pi}{2 \rho_1}}\right)^{1-\rho_1} e^{i \omega \left(s+te^{i \frac{\pi}{2 \rho_1}}\right)^{\rho_1}} \, dt \nonumber \\
												& \qquad \qquad + (1-\rho_1) \, e^{i \frac{\pi}{2 \rho_1}} \int_0^{+\infty} t^{N-1} \left(s+te^{i \frac{\pi}{2 \rho_1}}\right)^{-\rho_1} e^{i \omega \left(s+te^{i \frac{\pi}{2 \rho_1}}\right)^{\rho_1}} \, dt \bigg) \, . \nonumber
		\end{align}
		The boundary terms in \eqref{ipp} are zero; indeed one can see that the term at $0$ vanishes and one can use \eqref{estinf} once again to study the term at infinity. Then by similar arguments to the ones of the preceding case, we obtain
		\begin{align}
			\Big| \phi_{N}^{(1)}(s,\omega,\rho_1,1) \Big|	& \leqslant \frac{\omega^{-1}}{\rho_1 (N-2)!} \int_0^{+\infty} t^{N-2} \left| \left(s+te^{i \frac{\pi}{2 \rho_1}}\right)^{1-\rho_1} e^{i \omega \left(s+te^{i \frac{\pi}{2 \rho_1}}\right)^{\rho_1}} \right| \, dt \nonumber \\
															& \qquad + \frac{(\rho_1-1) \omega^{-1}}{\rho_1 (N-1)!} \int_0^{+\infty} t^{N-1} \left| \left(s+te^{i \frac{\pi}{2 \rho_1}}\right)^{-\rho_1} e^{i \omega \left(s+te^{i \frac{\pi}{2 \rho_1}}\right)^{\rho_1}} \right| \, dt \nonumber \\
															& \leqslant \frac{1}{\rho_1 (N-2)!} \, \omega^{-1} s^{1 - \rho_1} \int_0^{+\infty} t^{N-2} \, e^{- \omega t^{\rho_1}} \, dt \nonumber \\
															& \qquad + \frac{(\rho_1-1)}{\rho_1 (N-1)!} \, \omega^{-1} s^{-\rho_1} \int_0^{+\infty} t^{N-1} \, e^{- \omega t^{\rho_1}} \, dt \nonumber \\
															& = \frac{1}{\rho_1^2 (N-2)!} \, \Gamma\left( \frac{N-1}{\rho_1} \right) \, \omega^{-\left(1+ \frac{N-1}{\rho_1}\right)} s^{1-\rho_1} \nonumber \\
															& \qquad + \frac{(\rho_1 -1)}{\rho_1^2 (N-1)!} \, \Gamma \left( \frac{N}{\rho_1} \right) \, \omega^{-\left(1+\frac{N}{\rho_1}\right)} s^{-\rho_1} \nonumber \; .
		\end{align}
		Putting $b_{N,\rho_1} := \frac{1}{\rho_1^2 (N-2)!} \Gamma\Big( \frac{N-1}{\rho_1} \Big)$ and $c_{N,\rho_1} := \frac{(\rho_1 -1)}{\rho_1^2 (N-1)!} \Gamma \Big( \frac{N}{\rho_1} \Big)$, we can conclude this point.
	\end{itemize}
	A very similar work for $j=2$ provides the conclusion.
\end{proof}

\vspace{0.5cm}

Given a function which satisfies the two previous estimates, we furnish a new estimate in the second lemma, by employing the balance between blow-up and decay. Note that a technical argument requires $\rho \geqslant 2$.
\begin{lem2}\label{LEM2}
	Let $N \in \N \backslash \{0\}$, $\rho \geqslant 2$ and $f : (0,+\infty) \times (0,+\infty) \longrightarrow \R$ be a function which satisfies the two following inequalities:
	\begin{equation*}
		\forall \, s, \omega > 0 \qquad \left\{ \begin{array}{rl}
			& \displaystyle \big| f(s,\omega) \big| \leqslant a \, \omega^{-\frac{N}{\rho}} \; , \\
			& \displaystyle \big| f(s,\omega) \big| \leqslant b \, \omega^{-\left(1+\frac{N-1}{\rho}\right)} \, s^{1-\rho} \, + c \, \omega^{-\left(1+\frac{N}{\rho}\right)} \, s^{-\rho} \; , \\
		\end{array} \right.
	\end{equation*}
	where $a,b,c > 0$ are three constants.\\
	Fix $\gamma \in (0,1)$ and define $\delta := \rho^{-1} (\gamma+N) \in \Big(\frac{N}{\rho},\frac{1+N}{\rho} \Big)$. Then there exists a constant $L_{\gamma,\rho} > 0$ such that the following estimate holds:
	\begin{equation*}
		\forall \, s,\omega > 0  \qquad \big| f(s,\omega) \big| \leqslant L_{\gamma,\rho} \, s^{-\gamma} \, \omega^{-\delta} \; .
	\end{equation*}
	Furthermore $L_{\gamma,\rho} = a (K_{\rho})^{\gamma} > 0$ where $K_{\rho}$ is the unique positive solution of
	\begin{equation*}
		a K^{\rho} - b K - c = 0 \; .
	\end{equation*}	
\end{lem2}

\begin{proof}
	Let $g_1 , g_2 \, : \, (0,+\infty) \, \times \, (0,+\infty) \, \longrightarrow \R $ be defined as follows:
	\begin{equation*}
		g_1(s,\omega) := a \, \omega^{-\frac{N}{\rho}} \qquad , \qquad g_2(s,\omega):=  b \, \omega^{-\left(1+\frac{N-1}{\rho}\right)} \, s^{1-\rho} \, + c \, \omega^{-\left(1+\frac{N}{\rho}\right)} \, s^{-\rho} \; .
	\end{equation*}
	Now fix $\omega > 0$ and define the function $h_{\omega} : (0,+\infty) \longrightarrow \R$ by
	\begin{equation*}
		h_{\omega}(s) := s^{\rho} \big(g_1(s,\omega) - g_2(s,\omega) \big) = a \, \omega^{-\frac{N}{\rho}} \, s^{\rho} - b \, \omega^{-\left(1+\frac{N-1}{\rho}\right)} \, s - c \, \omega^{-\left(1+\frac{N}{\rho}\right)} \; .
	\end{equation*}
	One can ensure that the equation $h_{\omega}(s_{\omega})=0$ admits a unique positive solution $s_{\omega}$ given by $s_{\omega} := K_{\rho} \, \omega^{-\frac{1}{\rho}}$ where $K_{\rho}$ is the unique positive solution of the equation $a K^{\rho} - b K - c = 0$. So $g_1(.,\omega)$ and $g_2(.,\omega)$ intersect each other at the point $s_{\omega}$ and we have $g_1(.,\omega) \leqslant g_2(.,\omega)$ for $s \in (0,s_{\omega}]$ and $g_1(.,\omega) \geqslant g_2(.,\omega)$ otherwise. So we get more precise estimates:
	\begin{equation*}
	\left\{ \begin{array}{rl}
			& \displaystyle \forall \, s \in (0,s_{\omega}] \qquad \big| f(s,\omega) \big| \leqslant a \, \omega^{-\frac{N}{\rho}} = g_1(s,\omega) \; , \\
			& \displaystyle \forall \, s \in [s_{\omega},+\infty) \qquad \big| f(s,\omega) \big| \leqslant b \, \omega^{-\left(1+\frac{N-1}{\rho}\right)} \, s^{1-\rho} \, + c \, \omega^{-\left(1+\frac{N}{\rho}\right)} \, s^{-\rho} = g_2(s,\omega) \; . \\
		\end{array} \right.
	\end{equation*}
	Now let us build a function $g : (0,+\infty) \times (0,+\infty) \, \longrightarrow \R$ which is locally integrable with respect to the variable $s$ and which satisfies the following inequalities for any $\omega > 0$:
	\begin{equation} \label{int_cond}
		\left\{ \begin{array}{rl}
			& \displaystyle \forall \, s \in (0,s_{\omega}] \qquad \big| f(s,\omega) \big| \leqslant g_1(s,\omega) \leqslant g(s,\omega) \; , \\
			& \displaystyle \forall \, s \in [s_{\omega},+\infty) \qquad \big| f(s,\omega) \big| \leqslant g_2(s,\omega) \leqslant g(s,\omega) \; , \\
		\end{array} \right.
	\end{equation}
	Here we choose $g_{\gamma,\delta}(s,\omega) := L_{\gamma,\rho} \, s^{-\gamma} \omega^{-\delta}$, where $L_{\gamma,\rho}, \, \delta, \, \gamma > 0$ must be clarified. To this end, we require the following condition:
	\begin{equation*}
		\forall \, \omega > 0 \qquad g_{\gamma,\delta}(s_{\omega},\omega) = g_1(s_{\omega},\omega) = g_2(s_{\omega},\omega) \; .
	\end{equation*}
	This leads to
	\begin{equation*}
		g_{\gamma,\delta}\big(K_{\rho} \, \omega^{-\frac{1}{\rho}},\omega\big) = L_{\gamma,\rho} (K_{\rho})^{-\gamma} \omega^{\frac{\gamma}{\rho}-\delta} = a \, \omega^{-\frac{N}{\rho}} \; .
	\end{equation*}
	Since this equality holds for all $\omega > 0$, we obtain
	\begin{equation*}
		\left\{ \begin{array}{rl}
			& \displaystyle L_{\gamma,\rho} = a (K_{\rho})^{\gamma} \\
			& \displaystyle \frac{\gamma}{\rho} - \delta = -\frac{N}{\rho} \\
		\end{array} \right. \qquad \Longrightarrow \qquad
		\left\{ \begin{array}{rl}
			& \displaystyle L_{\gamma,\rho} = a (K_{\rho})^{\gamma} \\
			& \displaystyle \delta = \rho^{-1}(\gamma+N) \\
		\end{array} \right. \; .
	\end{equation*}
	We take $\gamma \in (0,1)$ to make $g_{\gamma, \delta}(.,\omega) : (0,+\infty) \longrightarrow \R$ locally integrable with respect to $s$; it follows $\delta = \rho^{-1} (\gamma+N) \in \Big(\frac{N}{\rho},\frac{1+N}{\rho} \Big)$. To conclude, we have to check inequalities \eqref{int_cond}.
	\begin{itemize}
		\item Take $s \leqslant s_{\omega}$. Then we have
		\begin{equation*}
			g_{\gamma,\delta}(s,\omega) = a (K_{\rho})^{\gamma} \omega^{-\delta} s^{-\gamma} \geqslant a \, \omega^{\frac{\gamma}{\rho} - \delta} = a \, \omega^{-\frac{N}{\rho}} = g_1(s,\omega) \; ,
		\end{equation*}
		since $\frac{\gamma}{\rho} - \delta = -\frac{N}{\rho}$.
		\item Choose $s \geqslant s_{\omega}$. We have to show that $g_2(s,\omega) \leqslant g_{\gamma,\delta}(s,\omega)$, that is to say
		\begin{equation} \label{toproove1}
			s^{\rho} \left( g_{\gamma,\delta}(s,\omega) - g_2(s,\omega) \right) = a (K_{\rho})^{\gamma} \omega^{-\delta} s^{\rho-\gamma} - b \, \omega^{-\left(1+\frac{N-1}{\rho}\right)} s - c \, \omega^{-\left(1+\frac{N}{\rho}\right)} \geqslant 0 \; .
		\end{equation}
		Define the function $k_{\omega} : (0,+\infty) \longrightarrow \R$ by $k_{\omega}(s) := s^{\rho}\big(g_{\gamma,\delta}(s,\omega) - g_2(s,\omega)\big)$, and differentiate it:
		\begin{equation*}
			(k_{\omega})'(s) = a (K_{\rho})^{\gamma} (\rho-\gamma) \, \omega^{-\delta} s^{\rho-\gamma -1} - b \, \omega^{-\left(1+\frac{N-1}{\rho}\right)} \; ,
		\end{equation*}
		Since $s > 0$ and $\rho \geqslant 2$, $(k_{\omega})'$ is an increasing function and vanishes at the point
		\begin{equation*}
			s_{\omega}' = \left( \frac{b}{a (K_{\rho})^{\gamma} (\rho - \gamma)} \right)^{\frac{1}{\rho - \gamma -1}} \omega^{-\frac{1}{\rho}} \; .
		\end{equation*}
		Now we want to show the inequality : $s_{\omega}' \leqslant s_{\omega}$. To this end, observe that
		\begin{equation*}
			0 \leqslant b K_{\rho} (\rho - \gamma - 1) + (\rho-\gamma) c \; ,
		\end{equation*}
		because $\rho \geqslant 2$. Furthermore since $K_{\rho}$ satisfies $a (K_{\rho})^{\rho} - b K_{\rho} - c = 0$, one gets
		\begin{equation*}
			\frac{b K_{\rho}}{\rho-\gamma} \leqslant b K_{\rho} + c = a (K_{\rho})^{\rho} \qquad \Longleftrightarrow \qquad \frac{b}{a (K_{\rho})^{\gamma} (\rho - \gamma)} \leqslant (K_{\rho})^{\rho - \gamma - 1} \; ,
		\end{equation*}
		and we obtain finally
		\begin{equation*}
			s_{\omega}' = \bigg(\frac{b}{a (K_{\rho})^{\gamma} (\rho - \gamma)}\bigg)^{\frac{1}{\rho-1-\gamma}} \omega^{-\frac{1}{\rho}} \leqslant K_{\rho} \, \omega^{-\frac{1}{\rho}} = s_{\omega} \; .
		\end{equation*}
		So for all $s \geqslant s_{\omega} \geqslant s_{\omega}'$, $k_{\omega}$ is an increasing function and thus we have
		\begin{equation*}
			k_{\omega}(s) \geqslant k_{\omega}(s_{\omega}) = (s_{\omega})^{\rho}\big(g_{\gamma,\delta}(s_{\omega},\omega) - g_2(s_{\omega},\omega)\big) = 0 \; .
		\end{equation*}
		Hence inequality \eqref{toproove1} is satisfied and $g_{\delta,\gamma}(s,\omega) - g_2(s,\omega) \geqslant 0$.
	\end{itemize}
\end{proof}	


\vspace{0.2cm}

We deduce the new estimates of the remainder from the two preceding results.

\begin{thm4} \label{NEW}
	Let $N \in \N \backslash \{0\}$, assume  $\mu_j = 1$ and $\rho_j \geqslant 2$ for a certain $j \in \{1,2\}$. Suppose that the functions $\psi : [p_1,p_2] \longrightarrow \R$ and $U : (p_1, p_2) \longrightarrow \C$ satisfy Assumption \emph{(P$_{\rho_1,\rho_2,N}$)} and Assumption \emph{(A$_{\mu_1,\mu_2,N}$)}, respectively. Then the statement of Theorem \ref{THM1} is still true and, for $\gamma \in (0,1)$ and $\delta := \rho_j^{-1}(\gamma + N) \in \big( \frac{N}{\rho_j}, \frac{N+1}{\rho_j} \big)$, we have more precise estimates for the remainder term $R_N^{(j)}(\omega)$:
	\begin{equation} \label{rest2}
		\left| R_N^{(j)}(\omega) \right| \leqslant L_{\gamma, \rho_j, N} \int_0^{s_j} s^{-\gamma} \left| \frac{d^N}{ds^N}(\nu_j k_j)(s) \right| \, ds \; \omega^{-\delta} \; ,
	\end{equation}
	for all $\omega > 0$, where $L_{\gamma, \rho_j, N} > 0$ is a constant given by Lemma \ref{LEM1} and Lemma \ref{LEM2}.
\end{thm4}

%
%
%
%
%
%

\begin{proof}
	We only have to check inequality \eqref{rest2} since the first four steps of the proof of Theorem \ref{THM1} remain valid with $\mu_j = 1$.\\
	Since Lemma \ref{LEM1} ensures that $ \phi_{N}^{(j)}(.,\omega,\rho_j,1)$ satisfies the assumptions of Lemma \ref{LEM2}, we get
	\begin{equation}
		\forall \, s \in (0,s_j] \quad  \forall \,\omega > 0 \qquad \left| \phi_{N}^{(j)}(s,\omega,\rho_j,1) \right| \leqslant L_{\gamma, \rho_j, N} \, s^{-\gamma} \, \omega^{-\delta} \; ,
	\end{equation}
	where $\gamma, \delta >0$ are defined as above and $L_{\gamma, \rho_j, N} > 0$ is given in Lemma \ref{LEM2}. Using the expression of the remainder term from Theorem \ref{THM1} and the preceding estimate leads to the conclusion, namely
	\begin{equation*}
		\begin{aligned}
			\left| R_{N}^{(j)}(\omega) \right|	& \leqslant \int_0^{s_j} \left| \, \phi_{N}^{(j)}(s,\omega,\rho_j, 1) \right| \, \left| \frac{d^N}{ds^N}(\nu_j k_j)(s) \right| ds \\
												& \leqslant L_{\gamma, \rho_j, N} \int_0^{s_j} s^{-\gamma} \left| \frac{d^N}{ds^N}(\nu_j k_j)(s) \right| ds \; \omega^{-\delta} \; .
		\end{aligned}
	\end{equation*}
	And we observe that the decay rate of the remainder term $R_N^{(j)}(\omega)$ with respect to $\omega$ is higher than the one of the highest term of the expansion $A_N^{(j)}(\omega)$.
\end{proof}


\section{Application to the free Schrödinger equation on the line: the influence of frequency domain singularities of the initial condition on dispersion}

\hspace{2.5ex} In this section, we are interested in the time-asymptotic behaviour of the solution of the free Schrödinger equation in one dimension. More precisely, we describe the motion and the localization of wave packets establishing asymptotic expansions. The influence on the decay of a singular frequency in the initial data is explored. We show that this singularity plays a key role in the dispersion and that the decay rates are affected.

We shall substantially use the results from the preceding section. Moreover the employed method is based on \cite[section 3]{fam}: we suppose that the initial data is in a frequency band. Hence the solution can be written as an oscillatory integral with respect to time with singular amplitude. Then an asymptotic expansion is obtained by applying the above mentioned stationary phase method. And finally we estimate uniformly the remainder.

In order to apply properly the results of Section 1, we shall distinguish the case where the stationary point of the phase belongs to the frequency band and where it does not. These considerations lead to two different expansions in two space-time regions: a cone generated by the frequency band and its outside.

Furthermore by applying directly the above stationary phase method, we show that the time decay rate is $t^{-\frac{\mu}{2}}$ on points moving in space-time with the critical velocity given by the singularity. In this case, we don't have to deal with the uniformity of the constant of the remainder since we establish an asymptotic expansion on a line.

Finally in the $L^2$ case, we use the previous expansion in the space-time cone to describe the asymptotic behaviour of the $L^2$-norm of the solution along the cross-section of this cone. The time decay rate $t^{-\frac{1}{2}}$ implies that this $L^2$-norm of the solution is asymptotically constant.

In this section, we explore only the case where the amplitude has a singularity. A similar work can be carried out if the amplitude is supposed non-singular by using the previous results with $\mu_1=\mu_2=1$.\\

Let us introduce the assumption concerning the initial data.\\

\noindent \textbf{Condition} ($C_{p_1,p_2,\mu}$): Fix $\mu \in (0,1)$ and let $p_1, p_2$ be two real numbers such that
	\begin{equation*}
		- \infty < p_1 < p_2 < +\infty \; .
	\end{equation*}
	$u_0$ satisfies condition ($C_{p_1,p_2,\mu}$) if and only if $\tf u_0 \equiv 0$ on $\R \setminus [p_1,p_2]$ and $U:=\tf u_0$ verifies Assumption (A$_{\mu,1,1}$) on $[p_1,p_2]$, with $\tf u_0 (p_2) = 0$.

\begin{rem4}
	\em \begin{enumerate}	\item If $u_0$ satisfies Condition $(C_{p_1,p_2,\mu})$ then the solution of the free Schrödinger equation $(S)$ is given by
\begin{equation*}
	u(t,x) = \frac{1}{2\pi} \int_{p_1}^{p_2} \tf u_0(p) \, e^{-itp^2 + ixp} \, dp \; ,
\end{equation*}
for all $t>0$ and $x \in \R$.

												\item The Fourier transform is formally defined by : $\displaystyle \tf u_0(p) = \int_{\R} u_0(x) \, e^{-i xp} \, dx \; .$
												\item Under Condition ($C_{p_1,p_2,\mu}$), $\tf u_0$ has a singularity of order $\mu-1$ at $p_1$ whereas the point $p_2$ is regular. For simplicity, we assume $\tf u_0 (p_2) = 0$ but a similar work can be carried out if $\tf u_0 (p_2) \neq 0$.
												\item $\tf u_0 \in \mathcal{C}^1\big((p_1,p_2]\big)$.
	\end{enumerate}
\end{rem4}

\vspace{0.3cm}

\begin{thm2} \label{THM2}
	Suppose that $u_0$ satisfies Condition ($C_{p_1,p_2,\mu}$). Fix $\delta \in \big(\max\{\mu,\frac{1}{2}\},1\big)$ and $\varepsilon_1, \varepsilon_2 > 0$ such that:
	\begin{equation*}
		p_1 + \varepsilon_1 < p_2 - \varepsilon_2 \; .
	\end{equation*}
	Then for all $(t,x) \in \mathfrak{C}_{\varepsilon_1,\varepsilon_2}(p_1,p_2)$ where $\mathfrak{C}_{\varepsilon_1,\varepsilon_2}(p_1,p_2)$ is the cone defined by
	\begin{equation*}
		p_1 + \varepsilon_1 \leqslant \frac{x}{2t} \leqslant p_2 - \varepsilon_2 \; ,
	\end{equation*}
	where $t > 0$, there exist $K_{\mu}(t,x,u_0), H(t,x,u_0) \in \C$ satisfying
	\begin{equation*}
		\big| u(t,x) - H(t,x,u_0) \, t^{-\frac{1}{2}} - K_{\mu}(t,x,u_0) \, t^{-\mu} \big| \leqslant c(u_0,\varepsilon_1, \varepsilon_2, \nu,\tilde{\nu}, \delta) \, \big(t^{-1} +  t^{-\delta} \big) \; ,
	\end{equation*}
	where $c(u_0,\varepsilon_1,\varepsilon_2,\nu,\tilde{\nu}, \delta) \geqslant 0$ is a constant independent on $t,x$. $\nu$ and $\tilde{\nu}$ are indexed families of smooth cut-off functions, introduced in the proof.
\end{thm2}

\begin{rem7} \label{REM4} \em See \eqref{FT1}, \eqref{FT2} and \eqref{R1} for explicit expressions of $K_{\mu}(t,x,u_0)$, $H(t,x,u_0)$ and $c(u_0,\varepsilon_1,\varepsilon_2,\nu,\tilde{\nu},\delta)$.
												
\end{rem7}
												
\begin{proof}
	We divide the proof in five steps, using a method based on \cite[section 3]{fam}.\\
	
	\noindent \textit{First step: Rewriting.} We factorize the phase function $p \longmapsto -tp^2 + xp$ by $t$, which gives
	\begin{equation*}
		\forall \, (t,x) \in [0,+\infty) \times \R \qquad u(t,x) = \frac{1}{2\pi} \int_{p_1}^{p_2} \tf u_0(p) e^{i t \Psi(p,t,x)} dp \; ,
	\end{equation*}
	where $\Psi(p,t,x) := - p^2 + \frac{x}{t} p = -\big( p - \frac{x}{2t} \big)^2 + \frac{x^2}{4 t^2}$. \footnote[1]{See Remark \ref{REM0}} \\
	
	\noindent \textit{Second step: Splitting of the integral.} We define the two functions $U$ and $\psi$ by
	\begin{equation*}
	\left\{ \begin{array}{rl}
			& \displaystyle \forall \, p \in (p_1,p_2] \qquad U(p) := \tf u_0(p) = (p-p_1)^{\mu-1} \tilde{u}(p) \; , \\
			& \displaystyle \forall \, p \in [p_1,p_2] \qquad \psi(p) := \Psi(p,t,x) \; .
	\end{array} \right.
	\end{equation*}
	Hence $\psi'$ can be written as follows:
	\begin{equation} \label{psi}
		\forall \, p \in [p_1,p_2] \qquad \psi'(p) = - 2 \left(p - \frac{x}{2 t}\right) = -2 \big(p-p_0(t,x) \big) \; ,
	\end{equation}
	where we put $p_0(t,x) := \frac{x}{2t}$. Henceforth, we will denote $p_0(t,x)$ by $p_0$. Thus $p_0$ is the only stationary point of the phase and it belongs to $[p_1 + \varepsilon_1, p_2 - \varepsilon_2]$ if and only if $(t,x) \in \mathfrak{C}_{\varepsilon_1,\varepsilon_2}(p_1,p_2)$. In such a situation, we write
	\begin{equation*}
		\begin{aligned}
			\forall \, (t,x) \in \mathfrak{C}_{\varepsilon_1,\varepsilon_2}(p_1,p_2) \qquad u(t,x)	& = \frac{1}{2\pi} \bigg( \int_{p_1}^{p_0} U(p) e^{i t \psi(p)} dp \; + \; \int_{p_0}^{p_2} U(p) e^{i t \psi(p)} dp \bigg) \\
		& =: \frac{1}{2\pi} \big( I^{(1)}(t,p_0) + I^{(2)}(t,p_0) \big) \; .
		\end{aligned}
	\end{equation*}
	
	\noindent \textit{Third step: Application of the stationary phase method.} The amplitude $U$ and the phase $\psi$ of $I^{(1)}(t,p_0)$ satisfy Assumption (A$_{\mu,1,1}$) and Assumption (P$_{1,2,N}$) for all $N \geqslant 1$ on $[p_1,p_0]$, respectively. This allows to apply Theorem \ref{THM1} to $I^{(1)}(t,p_0)$, where $t$ is the large parameter. But we must be careful with the cut-off function $\nu$ introduced in Theorem \ref{THM1}; effectively it depends on the point $p_0$ (and so on $t,x$) since the interval of integration depends also on $p_0$. Consequently we must consider an indexed family of smooth functions $\nu := \big\{ \nu_{p_0} \big\}_{p_0 \in [p_1 + \varepsilon_1, p_2 - \varepsilon_2]}$ defined as follows: choose a sufficiently small number $\eta_1$, for example $\eta_1 := \frac{\varepsilon_1}{3}$, independent on $p_0$. Then for any $p_0 \in [p_1 + \varepsilon_1, p_2 - \varepsilon_2]$, the smooth cut-off function $\nu_{p_0}$ is defined on $[p_1,p_0]$ like in Definition \ref{DEF1} and we extend it to the interval $[p_0,p_2]$ by $0$. In this situation, every $\nu_{p_0}$ is bounded by $1$ on $[p_1,p_2]$ and the graph compresses as $p_0$ tends to $p_1 + \varepsilon_1$; so we may suppose that $\big\| (\nu_{p_0})' \big\|_{L^{\infty}(p_1,p_2)}$ reaches its maximum when $p_0 = p_1 + \varepsilon_1$ and, without loss of generality, we may claim that
	\begin{equation} \label{nu}
		\forall \, p_0 \in [p_1+\varepsilon_1, p_2-\varepsilon_2] \qquad \big\| (\nu_{p_0})' \big\|_{L^{\infty}(p_1,p_2)} \leqslant \big\| (\nu_{p_1 + \varepsilon_1})' \big\|_{L^{\infty}(p_1,p_2)} =: M_{p_1, \varepsilon_1} \; .
	\end{equation}
	That is to say, $\nu$ is a bounded family in the space $\mathcal{C}^1 \big( [p_1,p_2] \big)$.\\ Applying Theorem \ref{THM1} provides
	\begin{equation*}
			I^{(1)}(t,p_0) = e^{i t \psi(p_1)} \Theta_1^{(1)}(1,\mu) k_1(0) \, t^{-\mu} + e^{it \psi(p_0)} \Theta_1^{(2)}(2,1) k_2(0) \, t^{-\frac{1}{2}} + R_1^{(1)}(t,p_0) + R_1^{(2)}(t,p_0) \; ,
	\end{equation*}
	for all $t > 0$, with $\nu_{1,p_0} := \nu_{p_0} \circ \varphi_1^{-1}$ and $\nu_{2,p_0} := (1-\nu_{p_0}) \circ \varphi_2^{-1}$.\\
	Let us compute the first terms. We employ the writing of $k_j$ given in \eqref{defk1} to obtain
	\begin{equation} \label{k0}
		k_1(0) = (\varphi_1^{-1})'(0)^{\mu} \, \tilde{u}(p_1) \qquad \text{and} \qquad k_2(0) = (\varphi_2^{-1})'(0) \, U(p_0) \; .
	\end{equation}
	We use now the expressions of $\varphi_j$ furnished in Definition \ref{DEF1}; note that the expression of $\varphi_2$ can be simplified:
	\begin{equation*}
		\varphi_1(p) = \psi(p) - \psi(p_1) \qquad \text{and} \qquad \varphi_2(p) = p_0 - p \; ,
	\end{equation*}
	for all $p \in [p_1, p_0-\eta_1]$ and $p \in [p_1 + \eta_1, p_0]$ respectively. It follows:
	\begin{equation} \label{deriv4}
		(\varphi_1^{-1})'(s) = \psi'\big(\varphi_1^{-1}(s)\big)^{-1} \qquad \text{and} \qquad (\varphi_2^{-1})'(s)= -1 \; ,
	\end{equation}
	for all $s \in [0,s_1]$ and $s \in [0,s_2]$ respectively. Therefore, using the explicit formulas for $\Theta_1^{(j)}(\rho_j,\mu_j)$ (see Theorem \ref{THM1}) and equalities \eqref{k0}, \eqref{deriv4}, we obtain
	\begin{equation*}
		\begin{aligned}
			I^{(1)}(t,p_0)	& = e^{i t \psi(p_1)} \, \Gamma(\mu) \, e^{i \frac{\pi \mu}{2}} \, \frac{\tilde{u}(p_1)}{\big(2 (p_0-p_1)\big)^{\mu}} \, t^{-\mu} + R_1^{(1)}(t,p_0) \\
							& \qquad  + e^{i t \psi(p_0)} \, \frac{\sqrt{\pi}}{2} e^{-i \frac{\pi}{4}} \, U(p_0) \, t^{-\frac{1}{2}} + R_1^{(2)}(t,p_0) \\
							& = \frac{\Gamma(\mu)}{2^{\mu}} \, e^{i \frac{\pi \mu}{2}} \, e^{i(-tp_1^2 + xp_1)} \, \tilde{u}(p_1) \, \Big(\frac{x}{2t} - p_1 \Big)^{-\mu} \,  t^{-\mu} + R_1^{(1)}(t,p_0) \\
							& \qquad + \frac{\sqrt{\pi}}{2} \, e^{-i \frac{\pi}{4}} \,  e^{i \frac{x^2}{4t}} \, \tilde{u}\Big( \frac{x}{2t} \Big) \, \Big(\frac{x}{2t} - p_1 \Big)^{\mu-1} \, t^{-\frac{1}{2}} + R_1^{(2)}(t,p_0) \; .
		\end{aligned}
	\end{equation*}
	
	\noindent \textit{Fourth step: uniform estimates for the remainder.} To finish the study of $I^{(1)}(t,p_0)$, let us estimate uniformly in $x$ the two remainders. Since the amplitude of $I^{(1)}(t,p_0)$ has a singularity at $p_1$, we can control $R_1^{(1)}(t,p_0)$ by using the remainder estimate of Theorem \ref{THM1}:
	\begin{equation} \label{unifb2}
		\left| R_1^{(1)}(t,p_0) \right| \leqslant \int_0^{s_1} s^{\mu-1} \big|(\nu_{1,p_0} k_1)'(s) \big| \, ds \, t^{-1}  \leqslant \frac{1}{\mu} \, (s_1)^{\mu} \, \big\|(\nu_{1,p_0} k_1)' \big\|_{L^{\infty}(0,s_1)} \, t^{-1} \; .
	\end{equation}
	Concerning the second remainder, we note that $p_0$ is not a singular point of the amplitude. In this case, we employ Theorem \ref{NEW}. Let us choose $\delta \in \big(\frac{1}{2}, 1 \big)$ and set $\gamma := 2\delta + 1 \in (0,1)$; then we obtain
	\begin{equation} \label{unifb1}
		\left| R_1^{(2)}(t,p_0) \right| \leqslant L_{\gamma, 2, 1} \int_0^{s_2} s^{-\gamma} \big|(\nu_{2,p_0} k_2)'(s) \big| \, ds \, t^{-\delta} \leqslant \frac{L_{\gamma, 2, 1}}{1-\gamma} \, (s_2)^{1-\gamma} \, \big\|(\nu_{2,p_0} k_2)' \big\|_{L^{\infty}(0,s_2)} \, t^{-\delta} \; .
	\end{equation}
	To bound explicitly the right-hand sides of \eqref{unifb2} and \eqref{unifb1}, one can see that
	\begin{equation*}
		 \big\|(\nu_{j,p_0} k_j)' \big\|_{L^{\infty}(0,s_j)} \leqslant \big\| (\varphi_j^{-1})' \big\|_{L^{\infty}(0,s_j)} \big\| (\nu_{p_0})' \big\|_{L^{\infty}(p_1,p_2)} \big\| k_j \big\|_{L^{\infty}(0,s_j)} + \big\| (k_j)' \big\|_{L^{\infty}(0,s_j)} \; ,
	\end{equation*}
	for $j=1,2$. Here, $\big\| (\nu_{p_0})' \big\|_{L^{\infty}(p_1,p_2)}$ is uniformly bounded by the constant $M_{p_1,\varepsilon_1}$ according to \eqref{nu}.	Further equality \eqref{psi} implies that $\psi'$ can be uniformly bounded from below by $2 \eta_1$ and from above by $2 (p_2-p_1)$ for all $p \in [p_1, p_0- \eta_1]$. Thanks to this observation, we obtain
	\begin{equation} \label{estphi}
		\forall \, s \in [0,s_1] \qquad \big(2 (p_2-p_1) \big)^{-1} \leqslant (\varphi_1^{-1})'(s) \leqslant (2 \eta_1)^{-1} \; .
	\end{equation}
	Moreover we recall that $s_1 = \varphi_1(p_0 - \eta_1)$; using the mean value Theorem, it follows
	\begin{equation*}
		s_1 = \psi(p_0 - \eta_1) - \psi(p_1) \leqslant 2 (p_2-p_1) (p_0 - \eta_1 - p_1 ) \leqslant 2 (p_2 - p_1)^2 \; .
	\end{equation*}
	Employing the expression of $k_1$ given in \eqref{defk1}, the following estimate is true for any $s \in [0,s_1]$,
	\begin{align}
		\big| k_1(s) \big|	& = \left| \left( \int_0^1 (\varphi_1^{-1})'(sy) dy \right)^{\mu-1} \tilde{u}\big(\varphi_1^{-1}(s)\big) \big( \varphi_1^{-1} \big)'(s) \right| \nonumber \\
							& \leqslant \big( 2 (p_2-p_1)\big)^{1-\mu} \| \tilde{u} \|_{L^{\infty}(p_1,p_2)} ( 2 \eta_1)^{-1} \nonumber \\
							& \leqslant \frac{(p_2 - p_1)^{1-\mu}}{2^{\mu}} \,  \eta_1^{-1} \,  \| \tilde{u} \|_{W^{1,\infty}(p_1,p_2)} \nonumber \; .
	\end{align}
	Now the product rule applied to $k_1$ on $[0,s_1]$ gives
	\begin{equation*}
		\begin{aligned}
			(k_1)'(s)	& = (\mu-1) \left(\int_0^1 y (\varphi_1^{-1})''(sy) dy \right) \left( \int_0^1 (\varphi_1^{-1})'(sy) dy \right)^{\mu-2} \tilde{u} \big( \varphi_1^{-1}(s) \big) (\varphi_1^{-1})'(s) \\
						& \qquad + \left( \int_0^1 (\varphi_1^{-1})'(sy) dy \right)^{\mu-1} \tilde{u}'\big( \varphi_1^{-1}(s) \big) (\varphi_1^{-1})'(s)^2 \\
						& \qquad \qquad + \left( \int_0^1 (\varphi_1^{-1})'(sy) dy \right)^{\mu-1} \tilde{u}\big( \varphi_1^{-1}(s) \big) (\varphi_1^{-1})''(s) \; .
		\end{aligned}
	\end{equation*}
	Using the relation $(\varphi_1^{-1})''(s) = - (\varphi_1)''\big(\varphi_1^{-1}(s)\big) (\varphi_1^{-1})'(s)^3 = 2 (\varphi_1^{-1})'(s)^3$ and estimates \eqref{estphi}, we obtain
	\begin{equation*}
		\big| (\varphi_1^{-1})''(s)	\big| = 2 \, \big| (\varphi_1^{-1})'(s)^3 \big| \leqslant 2^{-2} \, \eta_1^{-3} \; .
	\end{equation*}
	We employ inequalities \eqref{estphi} once again to bound $(k_1)'$ as follows:
	\begin{equation*}
		\begin{aligned}
			\big|(k_1)'(s)\big|	& \leqslant \frac{(1-\mu)}{2} \, 2^{-2} \, \eta_1^{-3} \big(2 (p_2-p_1)\big)^{2-\mu} \, \| \tilde{u}  \|_{L^{\infty}(p_1,p_2)} (2 \eta_1)^{-1} \\
								& \qquad + \big(2 (p_2-p_1)\big)^{1-\mu} \, \| \tilde{u}'  \|_{L^{\infty}(p_1,p_2)} (2 \eta_1)^{-2} \\
								& \qquad \qquad + \big(2 (p_2-p_1)\big)^{1-\mu} \, \| \tilde{u}  \|_{L^{\infty}(p_1,p_2)} 2^{-2} \eta_1^{-3} \\
								& \leqslant \frac{(p_2-p_1)^{1-\mu}}{2^{1+\mu}} \bigg(\frac{(1-\mu)}{2} \, \eta_1^{-4} (p_2-p_1)	 + \eta_1^{-2} + \eta_1^{-3} \bigg) \, \| \tilde{u} \|_{W^{1,\infty}(p_1,p_2)} \; .
		\end{aligned}
	\end{equation*}
	Putting everything together, we obtain
	\begin{equation*}
		\begin{aligned}
			\left| R_1^{(1)}(t,p_0) \right|	& \leqslant \frac{1}{\mu} \, \big(2 (p_2-p_1)^2\big)^{\mu} \Bigg( M_{p_1,\varepsilon_1} (2 \eta_1)^{-1} \frac{(p_2 - p_1)^{1-\mu}}{2^{\mu}} \, \eta_1^{-1} \, \| \tilde{u} \|_{W^{1,\infty}(p_1,p_2)} \\
											& \quad + \frac{(p_2-p_1)^{1-\mu}}{2^{1+\mu}} \bigg(\frac{(1-\mu)}{2} \, \eta_1^{-4} (p_2-p_1) + \eta_1^{-2} + \eta_1^{-3} \bigg) \| \tilde{u} \|_{W^{1,\infty}(p_1,p_2)} \Bigg) t^{-1} \\
											& =: c_1(u_0,\varepsilon_1, \nu) \, t^{-1}
		\end{aligned}
	\end{equation*}
	where
	\begin{equation} \label{c1}
		c_1(u_0,\varepsilon_1, \nu) := \frac{(p_2-p_1)^{1+\mu}}{2 \mu} \, \eta_1^{-2} \bigg( M_{p_1,\varepsilon_1} + \frac{(1-\mu)}{2} \, \eta_1^{-2} (p_2-p_1) + 1 + \eta_1^{-1} \bigg) \| \tilde{u} \|_{W^{1,\infty}(p_1,p_2)} \; .
	\end{equation}
	Let us study the second remainder. By the second relation of \eqref{deriv4}, we have
	\begin{equation}
		\min_{s \in [0,s_2]}\big| (\varphi_2^{-1})'(s) \big| = 1 = \big\| (\varphi_2^{-1})' \big\|_{L^{\infty}(0,s_2)} \; .
	\end{equation}
	Furthermore the definition of $s_2 = \varphi_2(p_1+\eta_1)$ implies
	\begin{equation*}
		s_2 = p_0 - p_1 - \eta_1 \leqslant p_2 - p_1 \; .
	\end{equation*}
	Now we note that for all $p \in [p_1 + \eta_1,p_0]$, we have $\big| U(p) \big| \leqslant \eta_1^{\mu-1} \| \tilde{u} \|_{L^{\infty}(p_1,p_2)}$; so we bound $k_2$ as follows
	\begin{equation*}
		\big|k_2(s) \big| = \left| U \big( \varphi_2^{-1}(s) \big) (\varphi_2^{-1})'(s) \right| \leqslant \eta_1^{\mu-1} \| \tilde{u} \|_{L^{\infty}(p_1,p_2)} \; .
	\end{equation*}
	for all $s \in [0,s_2]$. One can observe that $(\varphi_2^{-1})''(s) = 0$, so by the product rule, we have
	\begin{equation*}
		\forall \, s \in [0,s_2] \qquad (k_2)'(s) = U'\big( \varphi_2^{-1}(s) \big) (\varphi_2^{-1})'(s)^2 \; .
	\end{equation*}
	Since
	\begin{equation*}
		\begin{aligned}
			\big| U'(p) \big|	& \leqslant \big| (\mu-1) (p-p_1)^{\mu-2} \, \tilde{u}(p) \big| + \big| (p-p_1)^{\mu-1} \, \tilde{u}'(p) \big| \\
								& \leqslant (1-\mu) \, \eta_1^{\mu-2} \, \| \tilde{u} \|_{L^{\infty}(p_1,p_2)} + \eta_1^{\mu-1} \, \| \tilde{u}' \|_{L^{\infty}(p_1,p_2)} \; ,
		\end{aligned}
	\end{equation*}
	holds for all $p \in [p_1 + \eta_1,p_0]$, we obtain
	\begin{equation*}
		\big| (k_2)'(s) \big| \leqslant (1-\mu_1) \, \eta_1^{\mu-2} \, \| \tilde{u} \|_{L^{\infty}(p_1,p_2)} + \eta_1^{\mu-1} \, \| \tilde{u}' \|_{L^{\infty}(p_1,p_2)} \; ,
	\end{equation*}
	for all $s \in [0,s_2]$. Put everything together:
	\begin{equation*}
		\begin{aligned}
			\left| R_1^{(2)}(t,p_0) \right|	& \leqslant \frac{L_{\gamma,2,1}}{1-\gamma} \, (p_2-p_1)^{1-\gamma} \Big( M_{p_1,\varepsilon_1} \, \eta_1^{\mu-1} \, \| \tilde{u} \|_{L^{\infty}(p_1,p_2)} \\
												& \qquad + (1-\mu_1) \, \eta_1^{\mu-2} \, \| \tilde{u} \|_{L^{\infty}(p_1,p_2)} + \eta_1^{\mu-1} \, \| \tilde{u}' \|_{L^{\infty}(p_1,p_2)} \Big) \, t^{-\delta} \\
												& \leqslant c_2(u_0,\varepsilon_1, \nu,\delta) \, t^{-\delta} \; ,
		\end{aligned}
	\end{equation*}
	with
	\begin{equation} \label{c2}
		c_2(u_0,\varepsilon_1, \nu,\delta) := \frac{L_{\gamma,2,1}}{1-\gamma} \, (p_2-p_1)^{1-\gamma} \, \eta_1^{\mu-1} \Big( M_{p_1,\varepsilon_1} + (1-\mu_1) \eta_1^{-1} + 1 \Big) \| \tilde{u} \|_{W^{1,\infty}(p_1,p_2)} \; ;
	\end{equation}		
	the constant $c_2(u_0,\varepsilon_1,\nu,\delta)$ depends on $\delta$ since $\gamma = 2 \delta - 1$.\\
	
	\noindent \textit{Fifth step: Application of the stationary phase method to $I^{(2)}(t,p_0)$ and uniform estimates.} It remains to study the integral $I^{(2)}(t,p_0)$. First of all, remark that $\psi$ is decreasing on $[p_0,p_2]$; for example, we can use the substitution $p \longmapsto -p$ to make it increasing. Then we observe that the new amplitude and the new phase verify Assumption (A$_{1,1,1}$) and Assumption (P$_{1,2,N}$) (for $N \geqslant 1$) on $[-p_2,-p_0]$, respectively. Moreover we note that $U(p_2) = 0$ which implies that the first term related to $p_2$ vanishes. In consequences, we can estimate the remainder term related to $p_2$ by using Theorem \ref{THM1}, which is sufficient in this situation; the refinement of Theorem \ref{NEW} is not needed here. So by a similar work, we obtain the following expansion for $I^{(2)}(t,p_0)$:
	\begin{equation*}
		\begin{aligned}
			\Big| I^{(2)}(t,p_0) - \frac{\sqrt{\pi}}{2} \, e^{-i \frac{\pi}{4}} \,  e^{i \frac{x^2}{4t}} \, \tilde{u}\Big(\frac{x}{2t} \Big) \Big(\frac{x}{2t} - p_1\Big)^{\mu-1}	& \, t^{-\frac{1}{2}} \Big| \leqslant  \\
												& c_3(u_0,\varepsilon_1, \varepsilon_2, \tilde{\nu}) \, t^{-1} + c_4(u_0,\varepsilon_1, \varepsilon_2, \tilde{\nu}, \delta) \, t^{-\delta} \; ,
		\end{aligned}
	\end{equation*}
	where $c_3(u_0,\varepsilon_1, \varepsilon_2,\tilde{\nu}, \delta), c_4(u_0,\varepsilon_1, \varepsilon_2, \tilde{\nu}) \geqslant 0$ are defined as follows,
	\begin{align}
		& \label{c3} \bullet \; c_3(u_0,\varepsilon_1, \varepsilon_2, \tilde{\nu}) := \frac{(p_2-p_1)^2}{2} \, \eta_2^{-2} \varepsilon_1^{\mu-1} \Big( M_{p_2,\varepsilon_2} + (1-\mu) \varepsilon_1^{-1} + 1 + \eta_2^{-1} \Big)  \| \tilde{u} \|_{W^{1,\infty}(p_1,p_2)} \; , \\
		& \label{c4} \bullet \; c_4(u_0,\varepsilon_1, \varepsilon_2, \tilde{\nu},\delta) := \frac{L_{\gamma,2,1}}{1-\gamma} \, (p_2-p_1)^{1-\gamma} \varepsilon_1^{\mu-1} \Big( M_{p_2,\varepsilon_2} + (1-\mu) \varepsilon_1^{-1}+ 1 \Big) \| \tilde{u} \|_{W^{1,\infty}(p_1,p_2)} \; ,
	\end{align}
	where $\tilde{\nu} = \big\{ \tilde{\nu}_{-p_0} \big\}_{-p_0 \in [-p_2 + \varepsilon_2, -p_1-\varepsilon_1]}$ is an indexed family of smooth cut-off functions on $[-p_2,-p_1]$, defined in a similar way than $\nu$ in the third step. Hence $\tilde{\nu}$ is a bounded family in $\mathcal{C}^1([-p_2,-p_1])$ and we put $M_{p_2,\varepsilon_2} := \big\| (\tilde{\nu}_{-p_2+\varepsilon_2})' \big\|_{L^{\infty}(-p_2,-p_1)}$.\\
	
	Choosing $\delta > \max\big\{\frac{1}{2}, \mu\big\}$ implies that the decay rate of the remainder terms are higher than the decay rate of the first terms. And we conclude the proof by giving the expressions for
	\begin{align}
		& \label{FT1} \bullet \quad H(t,x,u_0) := \frac{1}{2 \sqrt{\pi}} \, e^{-i \frac{\pi}{4}} \,  e^{i \frac{x^2}{4t}} \, \tilde{u}\Big(\frac{x}{2t} \Big) \Big(\frac{x}{2t} - p_1\Big)^{\mu-1} \; , \\
		& \label{FT2} \bullet \quad K_{\mu}(t,x,u_0) := \frac{\Gamma(\mu)}{2^{\mu+1} \pi} \, e^{i \frac{\pi \mu}{2}} \, e^{i(-tp_1^2 + xp_1)} \, \tilde{u}(p_1) \, \Big(\frac{x}{2t} - p_1 \Big)^{-\mu} \; ,
	\end{align}
	and using \eqref{c1}, \eqref{c2}, \eqref{c3} and \eqref{c4}, we set
	\begin{align}
		\bullet \quad c(u_0,\varepsilon_1, \varepsilon_2,\nu,\tilde{\nu}, \delta) :=	& \frac{1}{2\pi} \big( c_1(u_0,\varepsilon_1,\nu) + c_2(u_0,\varepsilon_1,\nu,\delta) \nonumber \\
							& \label{R1}  \qquad + c_3(u_0,\varepsilon_1, \varepsilon_2,\tilde{\nu}) + c_4(u_0,\varepsilon_1, \varepsilon_2,\tilde{\nu},\delta)\big) \; . 
	\end{align}
	This ends the proof.
\end{proof}

\begin{REM0} \label{REM0}
	\em At this stage, the authors of \cite{fam} introduced the large parameter $\omega := \sqrt{t^2 + x^2}$, and replaced $t$ and $x$ by the bounded parameters $\tau := \frac{t}{\omega}$ and $\chi := \frac{x}{\omega}$. This led to a family of phase functions which was globally bounded in $\mathcal{C}^4$ with respect to $\tau$ and $\chi$. This was necessary for the application of \cite{hormander}. In our context, it is sufficient to control the phase functions in space-time cones. Indeed the explicitness of our remainder estimate shows that its coefficient depends only on the quotient $\frac{x}{t}$, which is bounded in these cones. It is not necessary to have the global boundedness with respect to $t$ and $x$ separately. Therefore we can use $t$ as a large parameter instead of $\omega$, which is conceptually simpler and clearer.
\end{REM0}

\vspace{0.3cm}

\begin{thm3}
	Suppose that $u_0$ satisfies Condition ($C_{p_1,p_2,\mu}$). Choose $\varepsilon > 0$ such that
	\begin{equation*}
		-\frac{1}{\varepsilon} < p_1 - \varepsilon \qquad \text{and} \qquad  p_2 + \varepsilon < \frac{1}{\varepsilon} \; .
	\end{equation*}
	Then for all $(t,x)$ that lies in the cone $\mathfrak{C}_{1,\varepsilon}^c(p_1,p_2)$ (resp. $\mathfrak{C}_{2,\varepsilon}^c(p_1,p_2)$) defined by
	\begin{equation*}
		 -\frac{1}{\varepsilon} \leqslant \frac{x}{2t} \leqslant p_1 - \varepsilon \qquad \Big( \text{resp.} \quad  p_2 + \varepsilon \leqslant \frac{x}{2t} \leqslant \frac{1}{\varepsilon} \Big) \; ,
	\end{equation*}
	where $t > 0$, there exists $K_{1,\mu}^c(t,x,u_0)$ (resp. $K_{2,\mu}^c(t,x,u_0)$) satisfying
	\begin{equation*}
		\big| u(t,x) - K_{j,\mu}^c(t,x,u_0) \, t^{-\mu} \big| \leqslant c_j^c(u_0,\varepsilon) \, t^{-1} \quad , \quad j=1,2 \; ,
	\end{equation*}
	where $c_1^c(u_0,\varepsilon) \geqslant 0$ (resp. $c_2^c(u_0,\varepsilon) \geqslant 0$) is a constant independent on $t,x$.
\end{thm3}

\begin{proof}
	Recall the expression of $u(t,x)$:
	\begin{equation*}
		u(t,x) = \frac{1}{2\pi} \int_{p_1}^{p_2} \tf u_0(p) \, e^{-i tp^2 + ixp} \, dp = \frac{1}{2\pi} \int_{p_1}^{p_2} U(p) \, e^{i t \psi(p)} \, dp \; ,
	\end{equation*}
	where the last equality is obtained by using the two first steps of the proof of Theorem \ref{THM2}.\\
	Here one can observe that there is no stationary point in the interval $[p_1,p_2]$; indeed the phase has a unique stationary point $p_0 = \frac{x}{2t}$, and if $(t,x) \in \mathfrak{C}_{j,\varepsilon}^c(p_1,p_2)$ then $p_0$ does not belong to $[p_1,p_2]$. So in the two cases $j=1,2$, we do not have to split the integral as in the second step of the preceding proof. In the case $j=1$, we need a substitution to make the phase increasing. Afterwards one can apply Theorem \ref{THM1}, where $t$ is the large parameter, and where $U$ and $\psi$ satisfy the assumptions (A$_{1,\mu,1}$) and (P$_{1,1,N}$) (for $N \geqslant 1$) on $[-p_2,-p_1]$ respectively in the case $j=1$, and (A$_{\mu,1,1}$) and (P$_{1,1,N}$) (for $N \geqslant 1$) on $[p_1,p_2]$ in the other case. As above, $U(p_2)=0$ implies that we can use Theorem \ref{THM1} to estimate the remainder related to $p_2$. Moreover, note that the cut-off function used in Theorem \ref{THM1} will not depend on $p_0$ in this situation (since $p_0 \notin [p_1,p_2]$); so we can arbitrarily choose a function $\nu$ satisfying the assumptions of Definition \ref{DEF1} and we define $M := \| \nu' \|_{L^{\infty}(p_1,p_2)}$, independent on $t,x$. And similar calculations than those of the proof of Theorem \ref{THM2} lead to the result.\\
	We give the expressions:
	\begin{equation*}
		\begin{aligned}
			& \bullet \; K_{j,\mu}^c(t,x,u_0) := \frac{\Gamma(\mu)}{2^{\mu+1} \pi} \, e^{(-1)^j i \frac{\pi \mu}{2}} \, e^{i(-tp_1^2 + xp_1)} \, \tilde{u}(p_1) \, \bigg( (-1)^j \Big(\frac{x}{2t} - p_1 \Big)\bigg)^{-\mu} \; , \\
			& \bullet \; c_1^c(u_0,\varepsilon) := \frac{1}{4 \pi} (\varepsilon^{-1} + p_2) \varepsilon^{-2} \bigg( \mu^{-1} (p_2-p_1)^{\mu} \Big( M + \frac{1-\mu}{2} (\varepsilon^{-1} + p_2) \varepsilon^{-2} + 1 + \varepsilon^{-1} \Big) \\
			& \qquad \qquad \qquad \quad + \eta^{\mu-1} (p_2-p_1) \Big( M + (1-\mu) \eta^{-1} + 1 + \varepsilon^{-1} \Big) \bigg) \| \tilde{u} \|_{W^{1,\infty}(p_1,p_2)} \; , \\
			& \bullet \; c_2^c(u_0,\varepsilon) := \frac{1}{4 \pi} (\varepsilon^{-1} - p_1) \varepsilon^{-2} \bigg( \mu^{-1} (p_2-p_1)^{\mu} \Big( M + \frac{1-\mu}{2} (\varepsilon^{-1}-p_1) \varepsilon^{-2} + 1 + \varepsilon^{-1} \Big) \\
			& \qquad \qquad \qquad \quad+ \eta^{\mu-1} (p_2 - p_1) \Big( M + (1-\mu) \eta^{-1} + 1 + \varepsilon^{-1} \Big) \bigg) \|\tilde{u}\|_{W^{1,\infty}(p_1,p_2)} \; , \\
		\end{aligned}
	\end{equation*}
where the fixed number $\eta \in \big(0, \frac{p_2-p_1}{2}\big)$ is related to the function $\nu$.
\end{proof}

\vspace{0.5cm}

In the following result, we establish an asymptotic expansion of the solution on the critical direction given by the singularity. Since the singularity of the amplitude and the stationary point of the phase coincide, a slow decay rate is obtained.

\begin{thm5}
	Suppose that $u_0$ satisfies Condition $(C_{p_1,p_2,\mu})$. Then for all $(t,x)$ such that
	\begin{equation*}
		\forall \, t > 0 \qquad x = 2 p_1 \, t \; ,
	\end{equation*}
	there exists $L_{\mu}(t,u_0) \in \C$ satisfying
	\begin{equation*}
		\left| u(t,x) - L_{\mu}(t,u_0) \, t^{-\frac{\mu}{2}} \right| \leqslant c(u_0) \, \left(t^{-1} + t^{-\frac{1}{2}}\right) \; .
	\end{equation*}
	The coefficient $L_{\mu}(t,u_0)$ and the constant $c(u_0) \geqslant 0$ are defined by
	\begin{align}
			& \label{Lmu} \bullet \quad L_{\mu}(t,u_0) := \frac{1}{2} \, \Gamma \left( \frac{\mu}{2} \right) e^{-i \frac{\pi \mu}{4}} e^{i t p_1^2} \, \tilde{u}(p_1) \; , \\
			& \bullet \quad c(u_0) := \frac{(p_2-p_1)^2}{2} \, \eta^{\mu-3} \left\| \tilde{u} \right\|_{W^{1,\infty}(p_1,p_2)} \left( \left\| \nu' \right\|_{L^{\infty}(p_1,p_2)} + (2-\mu) \, \eta^{-1} + 1 \right) \nonumber \\
			& \qquad \qquad \quad + \frac{\sqrt{\pi}}{2 \mu} \, (p_2-p_1)^{\mu} \left\| \tilde{u} \right\|_{W^{1,\infty}(p_1,p_2)} \left( \left\| \nu' \right\|_{L^{\infty}(p_1,p_2)} + 1 \right) \; , \nonumber
	\end{align}
	where $\nu$ is a smooth cut-off function coming from Theorem \ref{THM1}, and the fixed number $\eta \in \big(0, \frac{p_2-p_1}{2}\big)$ is related to this function.
\end{thm5}

\begin{proof}
	Simple application of Theorem \ref{THM1}.
\end{proof}

\begin{rem6}
	\em A similar work shows that the decay rate on the direction $x = 2 p_2 \, t$ is given by $t^{-\min \{\mu,\frac{1}{2}\}}$.
\end{rem6}

\vspace{0.5cm}

The next result is a consequence of Theorem \ref{THM2}. It permits to describe the time-asymptotic behaviour of the $L^2$-norm of the solution on the spatial cross-section of the cone $\mathfrak{C}_{\varepsilon_1,\varepsilon_2}(p_1,p_2)$, assuming $u_0 \in L^2(\R)$.

\begin{cor2} \label{COR2}
	Suppose that $u_0$ satisfies Condition $(C_{p_1,p_2,\mu})$ with $\mu \in \big(\frac{1}{2},1 \big)$. Choose $\varepsilon_1, \varepsilon_2 > 0$ such that $p_1 + \varepsilon_1 < p_2-\varepsilon_2$. Then there exists a constant $c(u_0, \varepsilon_1, \varepsilon_2) \geqslant 0$ such that for all $t \geqslant 1$,
	\begin{equation*}
		\left| \left\| u(t,.) \right\|_{L^2(I_t)} - \frac{1}{\sqrt{2 \pi}} \, \left\| \tf u_0 \right\|_{L^2(p_1+\varepsilon_1, p_2 - \varepsilon_2)} \right| \leqslant c(u_0, \varepsilon_1, \varepsilon_2) \, t^{\frac{1}{2}-\mu} \; ,
	\end{equation*}
	where
	\begin{equation*}
		I_t := \left[ \, 2 \, (p_1 + \varepsilon_1) \,  t, \, 2 \, (p_2 - \varepsilon_2) \, t \right] \; .
	\end{equation*}
\end{cor2}
\begin{proof}
	Firstly, we apply Theorem \ref{THM2} in the case $\mu \in \big( \frac{1}{2}, 1 \big)$ to obtain for all $(t,x) \in \mathfrak{C}_{\varepsilon_1, \varepsilon_2}(p_1,p_2)$ and $t \geqslant 1$,
	\begin{align}
		\left| u(t,x) - H(t,x,u_0) \, t^{-\frac{1}{2}} \right|	& \leqslant \left| K_{\mu}(t,x,u_0) \right| \, t^{-\mu} + c(u_0,\varepsilon_1, \varepsilon_2) \, \big( t^{-1} + t^{-\delta} \big) \nonumber \\
																& \leqslant \left( \left| K_{\mu}(t,x,u_0) \right| + c(u_0,\varepsilon_1, \varepsilon_2) \, \big(t^{-1+\mu} + t^{-\delta+\mu} \big) \right) \, t^{-\mu} \nonumber \\
																& \label{estcor1} \leqslant \tilde{c}(u_0, \varepsilon_1, \varepsilon_2) \, t^{-\mu} \; ,
	\end{align}
	with $\delta \in (\mu, 1)$; according to Theorem \ref{THM2}, the constant $c(u_0,\varepsilon_1, \varepsilon_2)$ depends also on the choice of the cut-off functions and $\delta$, but we do not write it here. Note that we used the boundedness of the coefficient $K_{\mu}(t,x,u_0)$ in the cone $\mathfrak{C}_{\varepsilon_1, \varepsilon_2}(p_1,p_2)$. Now we integrate the square of inequality \eqref{estcor1} on $I_t$ to obtain
	\begin{align*}
		\left\| u(t,.) - H(t,x,u_0) \, t^{-\frac{1}{2}} \right\|_{L^2(I_t)}^2	& = \int_{I_t} \left| u(t,x) - H(t,x,u_0) \, t^{-\frac{1}{2}} \right|^2 \, dx \\
				& \leqslant \tilde{c}(u_0, \varepsilon_1, \varepsilon_2)^2 \, t^{-2\mu} \left| I_t \right| \\
				& = c(u_0, \varepsilon_1, \varepsilon_2)^2 \, t^{1-2\mu} \; ,
	\end{align*}
	where we put $c(u_0, \varepsilon_1, \varepsilon_2)^2 := 2 \, (p_2 - p_1 - \varepsilon_1 - \varepsilon_2) \, \tilde{c}(u_0, \varepsilon_1, \varepsilon_2)^2$. It follows:
	\begin{align*}
		\left| \left\| u(t,.) \right\|_{L^2(I_t)} - \left\| H(t,.,u_0) \right\|_{L^2(I_t)} t^{-\frac{1}{2}} \right|	& \leqslant \left\| u(t,.) - H(t,x,u_0) \, t^{-\frac{1}{2}} \right\|_{L^2(I_t)} \\
						& \leqslant c(u_0, \varepsilon_1, \varepsilon_2) \, t^{\frac{1}{2}-\mu} \; .
	\end{align*}
	Moreover by recalling the expression of $H(t,x,u_0)$ given in \eqref{FT1}, we get
	\begin{align*}
		\left\| H(t,.,u_0) \right\|_{L^2(I_t)}^2	& = \frac{1}{4 \pi} \, \int_{I_t} \left| \tilde{u}\left(\frac{x}{2t} \right) \left( \frac{x}{2t} - p_1 \right)^{\mu-1} \right|^2 \, dx \\
												& = \frac{t}{2 \pi} \int_{p_1 + \varepsilon_1}^{p_2 - \varepsilon_2} \left| \tilde{u}(y) (y-p_1)^{\mu-1} \right|^2 \, dy \\
												& = \frac{t}{2 \pi} \left\| \tf u_0 \right\|_{L^2(p_1 + \varepsilon_1, p_2 - \varepsilon_2)}^2 \; .
	\end{align*}
	The proof is now complete.
\end{proof}

\section{Technical aspects of oscillation control}

\hspace{2.5ex} In this last section, we state and show several results, used in the proof of Theorem \ref{THM1} but not proved in the original paper.

Throughout this section, $\omega > 0$ will be considered as a fixed real number. Moreover $j$ will belong to the set $\{1,2\}$ and we shall prove the propositions in the case $j=1$; the proofs in the case $j=2$ are similar and require only appropriate changes of calculations.\\

The first result shows that a specific complex exponential with a large parameter can be estimated on a line in the complex plane. The choice of this line is strongly related to the phase function of the exponential.

\begin{prop2} \label{PROP2}
	Fix $s > 0$ and $\rho_j \geqslant 1$. Let $\Lambda^{(j)}(s)$ be the curve of the complex plane introduced in Definition \ref{DEF1}. Then we have
	\begin{equation*}
		\forall \, z = s + t e^{(-1)^{j+1} i \frac{\pi}{2 \rho_j}} \in \Lambda^{(j)}(s) \qquad \left| \, e^{(-1)^{j+1} i \omega z^{\rho_j}} \right| \leqslant e^{- \omega t^{\rho_j}} \; .
	\end{equation*}
\end{prop2}

\begin{proof}
	 First of all, fix $s,t \geqslant 0$ and note that the equality
	\begin{equation} \label{inteq}
		i \rho_1 \omega \int_0^s \left( \xi + t e^{i \frac{\pi}{2 \rho_1}} \right)^{\rho_1 -1} \, d\xi = i \omega z^{\rho_1} + \omega t^{\rho_1}
	\end{equation}
	holds for all $z = s + t e^{i \frac{\pi}{2 \rho_1}}$ by a simple calculation. Moreover one can see that
	\begin{equation*}
		\forall \, \xi \in [0,s] \qquad 0 \leqslant \text{Arg}\left(\xi + te^{i \frac{\pi}{2 \rho_1}}\right) \leqslant \frac{\pi}{2 \rho_1} \; ;
	\end{equation*}
	and since $\rho_1 \geqslant 1$, it follows:
	\begin{equation*}
		0 \leqslant \text{Arg}\left( \left(\xi + te^{i \frac{\pi}{2 \rho_1}} \right)^{\rho_1-1} \right) \leqslant \frac{\pi(\rho_1-1)}{2 \rho_1} \leqslant \frac{\pi}{2} \; .
	\end{equation*}
	In consequences the imaginary part of the complex number $\big(\xi + te^{i \frac{\pi}{2 \rho_1}}\big)^{\rho_1-1}$ is positive and so the real part of the right-hand side in \eqref{inteq} is negative. Hence we get
	\begin{equation*}
		\Re{\big( i \omega z^{\rho_1} + \omega t^{\rho_1} \big)} \leqslant 0 \quad \Longrightarrow \quad \big| \, e^{i \omega z^{\rho_1}} \big| \, e^{\omega t^{\rho_1}} = \big| \, e^{i \omega z^{\rho_1} + \omega t^{\rho_1}} \big| = e^{\Re\left({i \omega z^{\rho_1} + \omega t^{\rho_1}} \right)} \leqslant 1 \; ,
	\end{equation*}
	which yields the result in the case $j=1$. To treat the case $j=2$, we use the following equality
	\begin{equation*}
		-i \rho_2 \, \omega \int_0^s \left( \xi + t e^{-i \frac{\pi}{2 \rho_2}} \right)^{\rho_2 -1} \, d\xi = -i \omega z^{\rho_2} + \omega t^{\rho_2} \; ,
	\end{equation*}
	and we carry out a similar work. This ends the proof.
\end{proof}

\vspace{0.5cm}

The second proposition assures that the substitution used in Theorem \ref{THM1} is available.

\begin{prop3} \label{PROP3}
	Let $\psi : [p_1,p_2] \longrightarrow \R$ be a function which satisfies Assumption \emph{(P$_{\rho_1,\rho_2,N}$)}. Consider the function $\varphi_j : I_j \longrightarrow \R$ introduced in Definition \ref{DEF1}. Then $\varphi_j$ is a $\mathcal{C}^{N+1}$-diffeomorphism between $I_j$ and $[0,s_j]$.
\end{prop3}

\begin{proof}
	First of all, we check that $\varphi_1 \in \mathcal{C}^{N+1}(I_1)$. To do so, recall that $\psi'(p) = (p-p_1)^{\rho_1-1} \tilde{\psi}_2(p)$, where we put $\tilde{\psi}_2(p):=(p_2-p)^{\rho_2-1} \tilde{\psi}(p)$ which belongs to $\mathcal{C}^N(I_1)$. Applying Taylor's Theorem with the integral form of the remainder to $\psi'$, we obtain the useful integral representation of $\varphi_1$:
	\begin{equation*}
		\forall \, p \in I_1 := [p_1 , p_2 -\eta] \qquad \varphi_1(p) = (p-p_1) \bigg( \int_0^1 y^{\rho_1 -1} \, \tilde{\psi}_2 \big( y(p-p_1) + p_1 \big) \, dy \bigg)^{1/\rho_1} \; .
	\end{equation*}
	Fix $k \in \{1,\ldots,N\}$ and let us compute formally the $k$th derivative of the above expression by using the product rule:
	\begin{equation} \label{deriv}
		\begin{aligned}
			\varphi_1^{(k)}(p)	& = (p-p_1) \frac{d^k}{dp^k}\left[ \bigg( \int_0^1 y^{\rho_1 -1} \, \tilde{\psi}_2\big( y(p-p_1) + p_1 \big) \, dy \bigg)^{1/\rho_1} \right] \\
							& \qquad \qquad + k \, \frac{d^{k-1}}{dp^{k-1}} \left[ \bigg( \int_0^1 y^{\rho_1 -1} \, \tilde{\psi}_2\big( y(p-p_1) + p_1 \big) \, dy \bigg)^{1/\rho_1} \right] \; .
		\end{aligned}
	\end{equation}
	The positivity and the regularity of the function $\tilde{\psi}_2$ allow to differentiate $k$ times under the integral sign the function $\displaystyle J_1 : p \longmapsto \int_0^1 y^{\rho_1 -1} \tilde{\psi}_2\big( y(p-p_1) + p_1 \big) \, dy$. Hence the $k$ first derivatives of the composite function $\displaystyle p \longmapsto \left( \int_0^1 y^{\rho_1 -1} \tilde{\psi}_2\big( y(p-p_1) + p_1 \big) \, dy \right)^{1/\rho_1}$ exist and are continuous; in particular, the expression \eqref{deriv} is well-defined for all $p \in I_1$ and $\varphi_1^{(k)}$ is continuous. Concerning the $(N+1)$th derivative, we must be careful because we do not suppose that $\tilde{\psi}_2 \in \mathcal{C}^{N+1}(I_1)$. However we can formally apply the product rule once again for $k=N+1$:
	\begin{align}
		\varphi_1^{(N+1)}(p)	& = \label{deriv2} (p-p_1) \,  \frac{d^{N+1}}{dp^{N+1}}\left[ \bigg( \int_0^1 y^{\rho_1 -1} \, \tilde{\psi}_2\big( y(p-p_1) + p_1 \big) \, dy \bigg)^{1/\rho_1} \right] \\
								& \label{deriv3} \qquad +(N+1) \, \frac{d^{N}}{dp^{N}}\left[ \bigg( \int_0^1 y^{\rho_1 -1} \, \tilde{\psi}_2\big( y(p-p_1) + p_1 \big) \, dy \bigg)^{1/\rho_1} \right] \; .
	\end{align}
	Note that the term \eqref{deriv3} is well-defined by the previous work. So let us study \eqref{deriv2}. Firstly define $h_1 : s \longmapsto s^{\mu_1-1}$ and recall the expression of the function $J_1$ defined above; by applying Fa\`a di Bruno's Formula to $h_1 \circ J_1$, we obtain
	\begin{equation*}
		\frac{d^{N+1}}{dp^{N+1}} \Big( h_1 \circ J_1 \Big)(p) = \sum C_N \underbrace{h_1^{(m_1+\ldots+m_{N+1})}\big(J_1(p)\big)}_{(i)} \prod_{l=1}^{N+1} \underbrace{\Big( J_1^{(l)}(p) \Big)^{m_l}}_{(ii)}
	\end{equation*}
	where the sum is over all the $(N+1)$-tuples $(m_1, ..., m_{N+1})$ satisfying: $1m_1+2m_2+3m_3+\ldots+(N+1)m_{N+1}=N+1$. We note that the term $(i)$ is well-defined by the positivity of $J_1$; moreover by the previous study, the term $(ii)$ is well-defined and continuous for any $l \neq N+1$. So we have to study $\Big( J_1^{(N+1)}(p) \Big)^{m_{N+1}}$ where $m_{N+1} \leqslant 1$ by the above constraint. Since the case $m_{N+1} = 0$ is clear, suppose that $m_{N+1} = 1$. To differentiate $J_1^{(N)}$, we differentiate first $J_1$ $N$ times under the integral sign, then we make the substitution $y=\frac{s-p_1}{p-p_1}$ and finally we apply the fundamental Theorem of calculus:
	\begin{equation*}
		\begin{aligned}
			J_1^{(N+1)}(p)	& = \frac{d}{dp} \left[ \int_0^1 y^{N+\rho_1-1} \big(\tilde{\psi}_2 \big)^{(N)}\big((p-p_1)y+p_1\big) \, dy \right] \\
							& = \frac{d}{dp} \left[ \frac{1}{(p-p_1)^{\rho_1+N}} \int_{p_1}^p (s-p_1)^{\rho_1+N-1} \big(\tilde{\psi}_2 \big)^{(N)}(s) \, ds \right] \\
							& = \frac{-(\rho_1 +N)}{(p-p_1)^{\rho_1+N+1}} \int_{p_1}^p (s-p_1)^{\rho_1+N-1} \big(\tilde{\psi}_2 \big)^{(N)}(s) \, ds \\
							& \qquad + \frac{1}{(p-p_1)^{\rho_1+N}} (p-p_1)^{N+\rho_1 -1} \, \big(\tilde{\psi}_2 \big)^{(N)}(p) \\
							& = \frac{-(\rho_1 +N)}{(p-p_1)} \int_0^1 y^{\rho_1+N-1} \big(\tilde{\psi}_2 \big)^{(N)}\big(y(p-p_1)+p_1 \big) \, dy \\
							& \qquad  + \frac{1}{(p-p_1)} \, \big(\tilde{\psi}_2 \big)^{(N)}(p) \; .
		\end{aligned}
	\end{equation*}
	Multiplying this equality by $(p-p_1)$, we observe that $p \in I_1 \longmapsto (p-p_1) J_1^{(N+1)}(p)$ is well-defined and continuous, so is expression \eqref{deriv2}. Then $\varphi_1^{(N+1)}$ is continuous on $I_1$ and these considerations prove that $\varphi_1 \in \mathcal{C}^{N+1}(I_1)$.\\
	Furthermore, one remarks that
	\begin{equation*}
	\left\{ \begin{array}{rl}
			& \displaystyle \forall \, p \in I_1 \backslash \{ p_1 \} \qquad \varphi_1'(p) = \frac{1}{\rho_1} \, \psi'(p) \, \big(\psi(p) - \psi(p_1) \big)^{\frac{1}{\rho_1}-1} > 0 \; , \\
			& \displaystyle \varphi_1'(p_1) = \frac{1}{\rho_1^{\rho_1}} \, \tilde{\psi}_2 (p_1)^{\frac{1}{\rho_1}} > 0 \; ,
	\end{array} \right. \; ,
\end{equation*}
	so by the inverse function theorem, $\varphi_1 \, : \, I_1 \, \longrightarrow \, [0,s_1]$ is a $\mathcal{C}^{N+1}$-diffeomorphism.
\end{proof}

\vspace{0.5cm}

The regularity of the functions $k_j$ given in Definition \ref{DEF1} is proved in the next result.

\begin{prop4} \label{PROP4}
	Let $U : (p_1,p_2) \longrightarrow \C$ be a function which satisfies Assumption \emph{(A$_{\mu_1,\mu_2,N}$)}. Consider the function $k_j : (0,s_j] \longrightarrow \C$ introduced in Definition \ref{DEF1}. Then $k_j$ can be extended to the interval $[0,s_j]$ and $k_j \in \mathcal{C}^{N}\big([0,s_j] \big)$.
\end{prop4}

\begin{proof}
	Define $\tilde{u}_2(p) := (p_2-p)^{\mu_2-1} \tilde{u}(p)$ for any $p \in \, [p_1,p_2)$. We have by the definition of $k_1$,
	\begin{align}
		k_1(s)	& = \big( \varphi_1^{-1}(s) - \varphi_1^{-1}(0) \big)^{\mu_1 -1} \, \tilde{u}_2 \big(\varphi_1^{-1}(s) \big) \, s^{1-\mu_1} \, \big(\varphi_1^{-1}\big)'(s) \nonumber \\
				& = \bigg( \frac{\varphi_1^{-1}(s) - \varphi_1^{-1}(0)}{s} \bigg)^{\mu_1-1} \tilde{u}_2 \big(\varphi_1^{-1}(s) \big) \big(\varphi_1^{-1}\big)'(s) \nonumber \\
				& \label{defk1} = \bigg( \int_0^1 \big( \varphi_1^{-1} \big)'(sy) \, dy \bigg)^{\mu_1-1} \tilde{u}_2 \big(\varphi_1^{-1}(s) \big) \big(\varphi_1^{-1}\big)'(s)
	\end{align}
	with $s \in \, (0,s_1]$. This relation holds in fact for all $s \in [0,s_1]$ with $\displaystyle k_1(0) = \tilde{u}_2(p_1) \big(\varphi_1^{-1}\big)'(0)^{\mu_1}$. The regularity of $k_1$ comes from the regularity of $\tilde{u}$, the positivity of $(\varphi_1^{-1})'$ and Proposition \ref{PROP3}.
\end{proof}

\vspace{0.5cm}

The three last propositions are devoted to the successive primitives of the function $s \longmapsto s^{\mu_j-1} \, e^{(-1)^{j+1} \omega s^{\rho_j}}$. We shall use complex analysis to obtain primitives with integral representations; the path of integration is the line $\Lambda^{(j)}(s)$, introduced in Definition \ref{DEF1}, on which we are able to control the oscillations of the integrands (see Proposition \ref{PROP2}).

In the following result, we construct a primitive of a holomorphic function which is related to $s \longmapsto s^{\mu_j-1} \, e^{(-1)^{j+1} \omega s^{\rho_j}}$; this primitive is given by an integral on the line $\Lambda^{(j)}(s)$. To this end, we construct a sequence of primitives such that each one is given by an integral on a finite path, and the sequence of these paths tends to $\Lambda^{(j)}(s)$. Then we show that this sequence of functions converges uniformly on every compact set which implies by a theorem of Weierstrass that its limit is the desired primitive.

\begin{prop5} \label{PROP5}
	Let $s_j > 0$ and $l \in (0, 1)$. Define the domains $D_j \subset \C$ and $U \subset \C$ as follows:
	\begin{equation*}
		\begin{aligned}
			& \bullet \; D_j := \Big\{ v^* + t_ve^{(-1)^{j+1} i \frac{\pi}{2\rho_j}} \in \C \, \Big| \, v^* \in (0, s_j + l) \; , \; |t_v| < l \Big\} \\
			& \bullet \; U := \C \setminus \Big\{z \in \C \, \big| \, \Re(z) \leqslant 0 \; , \; \Im(z) =0 \Big\}
		\end{aligned}
	\end{equation*}
	Fix $\mu_j \in (0,1]$, $\rho_j \geqslant 1$ and $n \in \N$; let $F_{n,\omega}^{(j)}(.,.) : U \times \C \longrightarrow \C$ be the function defined by
	\begin{equation*}
		F_{n,\omega}^{(j)}(v,w) := \frac{(-1)^{n}}{n!} (v-w)^n v^{\mu_j - 1} e^{(-1)^{j+1} i \omega v^{\rho_j}} \; .
	\end{equation*}
	Then for every $w \in D_j$, $F_{n,\omega}^{(j)}(.,w)$ has a primitive $H_{n,\omega}^{(j)}(.,w)$ on $D_j$ given by
	\begin{equation*}
		H_{n,\omega}^{(j)}(v,w) := - \int_{\Lambda^{(j)}(v)} F_{n,\omega}^{(j)}(z,w) \, dz = \frac{(-1)^{n+1}}{n!} \int_{\Lambda^{(j)}(v)} (z-w)^n z^{\mu_j - 1} e^{(-1)^{j+1} i\omega z^{\rho_j}} dz \; ,
	\end{equation*}
	where $\Lambda^{(j)}(v)$ is the curve described by
	\begin{equation*}
		z \in \Lambda^{(j)}(v) \qquad \Longleftrightarrow \qquad \exists \, t \geqslant 0 \quad z= v + te^{(-1)^{j+1} i\frac{\pi}{2\rho_j}} \; .
	\end{equation*}
\end{prop5}

\begin{proof}
	Take $w \in D_1$ and $n \in \N$. Firstly, we have to ensure that the integral $H_{n,\omega}^{(1)}(v,w)$ is well-defined for every $v \in D_1$. Since $v \in D_1$, we can write $v= v^* + t_v e^{i\frac{\pi}{2\rho_1}}$ where $0 < v^* < s_1 + l$ and $-l < t_v < l$; we observe that
	\begin{equation*}
		-H_{n,\omega}^{(1)}(v,w) = \int_{\Lambda^{(1)}(v)} F_{n,\omega}^{(1)}(z,w) \, dz = \int_{\Lambda^{(1)}(v,v^*)} \ldots \; + \int_{\Lambda^{(1)}(v^*)} \ldots \; ,
	\end{equation*}	
	where $\Lambda^{(1)}(v,v^*)$ is the segment which starts from the point $v$ and goes to $v^*$, and $\Lambda^{(1)}(v^*)$ is given in the theorem. We furnish two parametrizations of these paths:
	\begin{equation*}
		\begin{aligned}
			& \forall \, t \in [-t_v,0] \qquad \lambda_{v,v^*}^{(1)}(t) := v^* - te^{i\frac{\pi}{2\rho_1}} \in \Lambda^{(1)}(v,v^*) \; , \\
			& \forall \, t \in [0,+\infty) \qquad \lambda_{v^*}^{(1)}(t) := v^* + te^{i\frac{\pi}{2\rho_1}} \in \Lambda^{(1)}(v^*) \; .
		\end{aligned}
	\end{equation*}
	We obtain
	\begin{align}
		\Big| F_{n,\omega}^{(1)} \big( \lambda_{v,v^*}^{(1)}(t), w \big) \Big|	&	\leqslant \frac{1}{n!} \left| v^* - t e^{i \frac{\pi}{2\rho_1}} - w \right|^n \left| v^* - t e^{i \frac{\pi}{2\rho_1}} \right|^{\mu_1-1} \bigg| e^{i \omega \left(v^* - t e^{i \frac{\pi}{2\rho_1}}\right)^{\rho_1}} \bigg| \nonumber \\
															&	\label{ineqrefv} \leqslant \frac{1}{n!} \sum_{k=0}^n \binom{n}{k} |v^* - w|^{n-k} \, (v^*)^{\mu_1-1} \, t^k e^{-\omega (-t)^{\rho_1}} \; ,
	\end{align}
	where \eqref{ineqrefv} comes from the binomial Theorem, Proposition \ref{PROP2} and the geometric observation:
	\begin{equation*}
		\left| v^* + t e^{i \frac{\pi}{2\rho_1}} \right| \geqslant v^* \; .
	\end{equation*}
	And a very similar calculation provides
	\begin{equation} \label{ineqref}
		\Big| F_{n,\omega}^{(1)} \big( \lambda_{v^*}^{(1)}(t), w \big) \Big|	\leqslant \frac{1}{n!} \sum_{k=0}^n \binom{n}{k} |v^* - w|^{n-k} \, (v^*)^{\mu_1-1} \, t^k e^{-\omega t^{\rho_1}} \; .
	\end{equation}
	Since \eqref{ineqrefv} and \eqref{ineqref} define integrable functions on $[-t_v,0]$ and $[0,+\infty)$ respectively, and since $\big| (\lambda_{v,v^*}^{(1)})'(t) \big| = \big| (\lambda_{v^*}^{(1)})'(t) \big| = 1$, the function $F_{n,\omega}^{(1)}(., w)$ is integrable on the paths $\Lambda^{(1)}(v,v^*)$ and $\Lambda^{(1)}(v^*)$ and hence, $H_{n,\omega}^{(1)}(v,w)$ is well-defined.\\
	Now we want to show that $H_{n,\omega}^{(1)}(.,w) \, : \, D_1 \, \longrightarrow \, \C$ is a primitive of $F_{n,\omega}^{(1)}(.,w)$ on $D_1$. To this end, we show that it is a uniform limit on all compact subsets of $D_1$ of a sequence of functions $\big(H_{m,n,\omega}^{(1)}(.,w)\big)_{m\geqslant 1}$ which are primitives of $F_{n,\omega}^{(1)}(.,w)$ on $D_1$. Here we build this sequence as follows: fix an arbitrary point $v_0 \in (0,s_1 + l)$ and define the following sequence of complex numbers:
	\begin{equation*}
		\forall \, m \in \N \backslash \{0\} \qquad v_m := v_0 + m e^{i\frac{\pi}{2 \rho_1}} \; .
	\end{equation*}	
	Let $m \in \N \backslash \{0\}$, let $v = v^* + t_v e^{i\frac{\pi}{2 \rho_1}} \in D_1$ and let $\Lambda_m(v)$ be the path which is composed of the segment that starts from the point $v$ and goes to the point $v^* + me^{i\frac{\pi}{2\rho_1}}$ and the horizontal segment that joins the points $v^* + me^{i\frac{\pi}{2\rho_1}}$ and $v_m$. We can now define the sequence of functions $\big(H_{m,n,\omega}^{(1)}(.,w) : D_1 \longrightarrow \C\big)_{m\geqslant1}$ as follows:
	\begin{equation*}
		H_{m,n,\omega}^{(1)}(v,w) := -\int_{\Lambda_m(v)} F_{n,\omega}^{(1)}(z,w) \,  dz \; .
	\end{equation*}
	It is clear that $F_{n,\omega}^{(1)}(.,w)$ is holomorphic on $U$, which is simply connected, and for any $v \in D_1$, $\Lambda_m(v)$ is included in $U$. The Cauchy integral Theorem affirms that each $H_{m,n,\omega}^{(1)}(.,w) : D_1 \longrightarrow \C$ is a primitive of the function $F_{n,\omega}^{(1)}(.,w)$.\\ Now, we have to prove that this sequence converges to $H_{n,\omega}^{(1)}(.,w)$ uniformly on any compact subset $K$ of $D_1$. Let $K \subset D_1$ be compact and for every $v \in K$, we write:
	\begin{equation} \label{weier}
		H_{m,n,\omega}^{(1)}(v,w) - H_{n,\omega}^{(1)}(v,w) = \int_{\Lambda_m^{c,1}(v)} F_{n,\omega}^{(1)}(z,w) \, dz + \int_{\Lambda_m^{c,2}(v)} F_{n,\omega}^{(1)}(z,w) \, dz \; ,
	\end{equation}
	where $\Lambda_m^{c,1}(v)$ is the horizontal segment which starts from $v_m$ and goes to $v^* + m e^{i\frac{\pi}{2 \rho_1}}$, $\Lambda_m^{c,2}(v)$ is the half-line with angle $\frac{\pi}{2 \rho_1}$ that starts from $v^* + m e^{i\frac{\pi}{2 \rho_1}}$ and goes to infinity. Let $\lambda_m^{c,1} : [0,|v_0 - v^*|] \longrightarrow \C$ and $\lambda_m^{c,2} : [0,+\infty) \longrightarrow \C$ be two parametrizations of $\Lambda_{1,m}^c(v)$ and $\Lambda_{2,m}^c(v)$ respectively and defined by
	\begin{equation*}
		\begin{aligned}
			& \forall \, t \in [0,|v_0-v^*|] \qquad \lambda_m^{c,1}(t) := \pm t + v_0 + m e^{i\frac{\pi}{2\rho_1}} \in \Lambda_m^{c,1}(v) \; , \\
			& \forall \, t \in [0,+\infty) \qquad \lambda_m^{c,2}(t) := v^* + (t+m) e^{i\frac{\pi}{2\rho_1}} \in \Lambda_m^{c,2}(v) \; .
		\end{aligned}
	\end{equation*}
	Then we have the following estimates:
	\begin{align}
		\Big| F_{n,\omega}^{(1)}(\lambda_m^{c,1}(t),w) \Big|	& \label{inequa5} \leqslant \frac{1}{n!} \sum_{k=0}^n \binom{n}{k} | v_0 -w |^{n-k} \left|\pm t + m e^{i\frac{\pi}{2 \rho_1}}\right|^k m^{\mu_1 -1} e^{-\omega m^{\rho_1}} \\
																& \label{inequa6} \leqslant \frac{1}{n!} \sum_{k=0}^n \binom{n}{k} | v_0 -w |^{n-k} \big(C_1(K) + m \big)^k m^{\mu_1 -1} e^{-\omega m^{\rho_1}}
	\end{align}
	\begin{itemize}
		\item \eqref{inequa5}: use the binomial Theorem, Proposition \ref{PROP2} and $\big| \lambda_m^{c,1}(t) \big| \geqslant m$;
		\item \eqref{inequa6}: employing the compactness of $K$, we have $0 \leqslant t \leqslant |v_0-v^*| \leqslant C_1(K)$, for a certain constant $C_1(K) > 0$;
	\end{itemize}
	Parametrizing the integral gives
	\begin{equation*}
		\begin{aligned}
			\bigg| \int_{\Lambda_m^{c,1}(v)} F_{n,\omega}^{(1)}(z,w) \, dz \bigg|	& \leqslant \int_0^{|v_0-v^*|} \frac{1}{n!} \sum_{k=0}^n \binom{n}{k} | v_0 -w |^{n-k} \big(C_1(K) + m \big)^k m^{\mu_1 -1} e^{-\omega m^{\rho_1}} \, dt \\
		&\leqslant \frac{1}{n!} \sum_{k=0}^n \binom{n}{k} | v_0 -w |^{n-k} \big(C_1(K) + m \big)^k m^{\mu_1 -1} e^{-\omega m^{\rho_1}} C_1(K) \\
		& \longrightarrow 0 \quad , \quad m \longrightarrow +\infty \; ,
		\end{aligned}
	\end{equation*}
	where we used the fact that $|v_0-v^*| \leqslant C_1(K)$ one more time; here, the convergence is uniform with respect to $v$. Furthermore,
	\begin{align}
		\Big| F_{n,\omega}^{(1)}(\lambda_m^{c,2}(t),w) \Big|	& \label{inequa1} \leqslant  \frac{1}{n!} \sum_{k=0}^n \binom{n}{k} |v^* - w|^{n-k} |v^*|^{\mu_1-1} (t+m)^k e^{-\omega (t+m)^{\rho_1}} \\
																& \label{inequa2} \leqslant \frac{C_{2,w}(K)}{n!} \, \sum_{k=0}^n \binom{n}{k} (t+m)^k e^{-\omega (t+m)^{\rho_1}} \\
																& \label{inequa3} \leqslant \frac{C_{2,w}(K)}{n!} \, \sum_{k=0}^n \binom{n}{k} m^k e^{-\omega m^{\rho_1}} (1+t)^k e^{-\omega t^{\rho_1}} \\
																& \label{inequa4} \leqslant \frac{C_{2,w}(K) \, M_{\omega}}{n!} \, \sum_{k=0}^n \binom{n}{k} (1+t)^k e^{-\omega t^{\rho_1}}
	\end{align}
	\begin{itemize}
		\item \eqref{inequa1}: use the binomial Theorem, Proposition \ref{PROP2} and $v^* \leqslant | \lambda_m^{c,2}(t) |$ ;
		\item \eqref{inequa2}: use the compactness of $K$ and $v \in K$ ;
		\item \eqref{inequa3}: $(m+t)^k \leqslant m^k (1+t)^k$ and $e^{-\omega(t+m)^{\rho_1}} \leqslant e^{-\omega m^{\rho_1}} e^{-\omega t^{\rho_1}}$ ;
		\item \eqref{inequa4}: use the boundedness of the sequences $(m^k e^{-\omega m^{\rho_1}})_{m\geqslant1}$ for $k=0,\ldots,n$.
	\end{itemize}
	We remark that \eqref{inequa3} tends to $0$ as $m$ tends to infinity for all $t \geqslant 0$ and \eqref{inequa4} gives an integrable function independent on $m$. So by the dominated convergence Theorem,
	\begin{equation*}
		\begin{aligned}
			\bigg| \int_{\Lambda_m^{c,2}(v)} F_{n,\omega}^{(1)}(z,w) \, dz \bigg|	& \leqslant \int_0^{+\infty} \Big| F_{n,\omega}^{(1)}(\lambda_m^{c,2}(t),w) \Big| \, dt \\
		& \leqslant \int_0^{+\infty} \frac{C_{2,w}(K)}{n!} \, \sum_{k=0}^n \binom{n}{k} m^k e^{-\omega m^{\rho_1}} (1+t)^k e^{-\omega t^{\rho_1}} \, dt \\
																					& \longrightarrow 0 \quad , \quad m \longrightarrow + \infty \; ,
		\end{aligned}
	\end{equation*}
	and the convergence is uniform with respect to $v$ since the last term is independent on $v$. Finally, these considerations and a theorem of Weierstrass imply that $H_{n,\omega}^{(1)}(.,w): D_1 \longrightarrow \C$ is a primitive of $F_{n,\omega}^{(1)}(.,w)$ on $D_1$.
\end{proof}

\vspace{0.5cm}

In the next step, we compute the successive primitives of the function $F_{n,\omega}^{(j)}(.,.)$ restricted to $\big\{(u,u) \, \big| \, u \in D_j \big\} \subset \C \times \C$. To this end, we use complex analysis in several variables and the preceding result.

\begin{prop6} \label{PROP6}
	Let $n \in \N \backslash \{0\}$. With the notations of Proposition \ref{PROP5}, define the function $h : \C \longrightarrow \C \times \C$ by
	\begin{equation*}
		h(u) := (u,u) \; ,
	\end{equation*}
	and let $H_{n,\omega}^{(j)}(.,.) : D_j \times D_j \longrightarrow \C$ be the function defined in Proposition \ref{PROP5}. Then the composite function $H_{n,\omega}^{(j)}(.,.) \circ h$ is holomorphic on $D_j$ and its derivative is given by
	\begin{equation*}
		\forall \, u \in D_j \qquad \frac{d}{d u} \big(H_{n,\omega}^{(j)}\circ h \big)(u) = \big(H_{n-1,\omega}^{(j)} \circ h \big)(u) \; .
	\end{equation*}
\end{prop6}

\begin{proof}
	The aim of the proof is to differentiate the composite function. For this purpose, we must ensure that this function is holomorphic with respect to each variable.\\
	Fix $n \in \N \backslash \{0\}$. We remark that each component of $h$ is holomorphic on $\C$, so is $h$ on $\C \times \C$. Moreover for any fixed $w \in  D_1$, $H_{n,\omega}^{(1)}(.,w) : D_1 \longrightarrow \C$ is a primitive of $F_{n,\omega}^{(1)}(.,w)$ on $D_1$ by Proposition \ref{PROP5}, so it is holomorphic. Now let us show that $H_{n,\omega}^{(1)}(v,.) : D_1 \longrightarrow \C$ belongs to $\mathcal{C}^1(D_1)$ and satisfies the Cauchy-Riemann Equations for fixed $v \in  D_1$. To do so we employ the holomorphy of $F_{n,\omega}^{(1)}(v,.) : \C \longrightarrow \C$ which provides the following relations:
	\begin{equation} \label{holomorphy}
		\forall \, w = x + iy \in \C \qquad  \frac{\partial}{\partial w} \big( F_{n,\omega}^{(1)} \big) (v,w) = \frac{\partial}{\partial x} \big( F_{n,\omega}^{(1)} \big) (v,w) = -i \frac{\partial}{\partial y} \big( F_{n,\omega}^{(1)} \big) (v,w) \; .
	\end{equation}
	And by a quick calculation, we get
	\begin{equation} \label{form1}
		\frac{\partial}{\partial w} \big( F_{n,\omega}^{(1)} \big) (v,w) = \frac{(-1)^{n-1}}{(n-1)!} \, (v-w)^{n-1} v^{\mu_1-1} e^{i \omega v^{\rho_1}} = F_{n-1,\omega}^{(1)}(v,w) \; .
	\end{equation}
	Furthermore, one can bound $F_{n-1,\omega}^{(1)}(.,w)$ on each path $\Lambda^{(1)}(v,v^*)$ and $\Lambda^{(1)}(v^*)$ by integrable functions independent on $w$. Indeed \eqref{ineqrefv} and \eqref{ineqref} show that $F_{n-1,\omega}^{(1)}(.,w)$ is bounded by integrable functions on $\Lambda^{(1)}(v,v^*)$ and $\Lambda^{(1)}(v^*)$, and the boundedness of $D_1$ allows to control $|w|$ by a constant in \eqref{ineqrefv} and \eqref{ineqref}. So we obtain the ability to differentiate under the integral sign which yields the following equalities:
	\begin{align}
		-\frac{\partial}{\partial x} \big( H_{n,\omega}^{(1)} \big) (v,w)	& = \frac{\partial}{\partial x} \left[ \int_{\Lambda^{(1)}(v,v^*)} F_{n,\omega}^{(1)}(z,w) \, dz \right] + \frac{\partial}{\partial x} \left[ \int_{\Lambda^{(1)}(v^*)} F_{n,\omega}^{(1)}(z,w) \, dz \right] \nonumber \\
												& = \label{equa1} \int_{\Lambda^{(1)}(v,v^*)} \frac{\partial}{\partial x} \big( F_{n,\omega}^{(1)} \big) (z,w) \, dz \: + \int_{\Lambda^{(1)}(v^*)} \frac{\partial}{\partial x} \big( F_{n,\omega}^{(1)} \big) (z,w) \, dz \\
												& = \label{equa2} \int_{\Lambda^{(1)}(v,v^*)} \frac{\partial}{\partial w} \big( F_{n,\omega}^{(1)} \big) (z,w) \, dz \: + \int_{\Lambda^{(1)}(v^*)} \frac{\partial}{\partial w} \big( F_{n,\omega}^{(1)} \big) (z,w) \, dz \\
												& = \label{equa3} \int_{\Lambda^{(1)}(v,v^*)} F_{n-1,\omega}^{(1)}(z,w) \, dz \: + \int_{\Lambda^{(1)}(v^*)} F_{n-1,\omega}^{(1)} (z,w) \, dz \\
												& = \int_{\Lambda^{(1)}(v)} F_{n-1,\omega}^{(1)}(z,w) \, dz \nonumber \\
												& = -H_{n-1,\omega}^{(1)}(v,w) \nonumber												
	\end{align}
	\begin{itemize}
		\item \eqref{equa1}: application of the Theorem of differentiation under the integral sign ;
		\item \eqref{equa2}: holomorphy of the function $F_{n,\omega}^{(1)}(v,.)$ by using equalities \eqref{holomorphy} ;
		\item \eqref{equa3}: relation \eqref{form1} .
	\end{itemize}
	In a similar way, we obtain
	\begin{equation*}
		-i \frac{\partial}{\partial y} \big( H_{n,\omega}^{(1)} \big) (v,w) = H_{n-1,\omega}^{(1)}(v,w) \; .
	\end{equation*}
	Then the Cauchy-Riemann Equations are satisified and $\displaystyle \frac{\partial}{\partial x} \big( H_{n,\omega}^{(1)} \big) (v,.)$ and $\displaystyle \frac{\partial}{\partial y} \big( H_{n,\omega}^{(1)} \big)(v,.)$ are continuous on $D_1$ by the continuity of $F_{n-1,\omega}^{(1)}(z,.) : \C \longrightarrow \C$. So $H_{n,\omega}^{(1)}(v,.) : D_1 \longrightarrow \C$ is holomorphic, with
	\begin{equation*}
		\frac{\partial}{\partial w} \big( H_{n,\omega}^{(1)} \big) (v,w) = \frac{\partial}{\partial x} \big( H_{n,\omega}^{(1)} \big) (v,w) = -i \frac{\partial}{\partial y} \big( H_{n,\omega}^{(1)} \big) (v,w) = H_{n-1,\omega}^{(1)}(v,w) \; .
	\end{equation*}
	Finally the composite function $H_{n,\omega}^{(1)} \circ h$ is holomorphic on $D_1 \times D_1$ and we have the formula
	\begin{align*}
		\frac{d}{d u} \big(H_{n,\omega}^{(1)} \circ h \big)(u)	& = \bigg( \frac{\partial}{\partial v } \big( H_{n,\omega}^{(1)} \big) \big( h(u) \big) \quad \frac{\partial}{\partial w} \big( H_{n,\omega}^{(1)} \big) \big(h(u) \big) \bigg) \left( \begin{array}{c}
						1 \\
						1 \\
						\end{array}
						\right) \\
																& =  \frac{\partial}{\partial v} \big( H_{n,\omega}^{(1)} \big) (u,u) +  \frac{\partial}{\partial w} \big( H_{n,\omega}^{(1)} \big) (u,u) \; ;
	\end{align*}
	And a short computation shows that $\displaystyle \frac{\partial}{\partial v} \big( H_{n,\omega}^{(1)} \big) (u,u) = F_{n,w}^{(1)}(u,u) = 0$, so
	\begin{equation*}
		\forall \, v \in D_1 \quad \frac{d}{d u} \big(H_{n,\omega}^{(1)}\circ h \big)(u) = \big(H_{n-1,\omega}^{(1)} \circ h \big)(u) = \frac{(-1)^n}{(n-1)!} \int_{\Lambda^{(1)}(u)} (z-u)^{n-1} z^{\mu_1-1} e^{i \omega z^{\rho_1}} dz \; .
	\end{equation*}
\end{proof}

\vspace{0.3cm}

From the two previous propositions, we deduce the final corollary.

\begin{cor1} \label{COR1}
	Fix $s_j > 0$, $\rho_j \geqslant 1$ and $\mu_j \in (0,1]$. For any $\omega >0$, the sequence of functions $\big( \phi_{n}^{(j)}(.,\omega,\rho_j,\mu_j) : \, (0,s_j] \longrightarrow \C \big)_{n \geqslant 1}$ defined in Theorem \ref{THM1} satisfies the recursive relation:
	\begin{equation*}
		\forall \, s \in \; (0,s_j] \qquad \left\{	\begin{array}{rl}
														& \displaystyle \frac{\partial}{\partial s} \big( \phi_{n+1}^{(j)} \big) (s,\omega,\rho_j,\mu_j) = \phi_{n}^{(j)}(s, \omega,\rho_j,\mu_j) \qquad \forall \, n \geqslant 1 \; , \\
														& \displaystyle \frac{\partial}{\partial s} \big( \phi_{1}^{(j)} \big)(s,\omega,\rho_j,\mu_j) = s^{\mu_j-1} e^{(-1)^{j+1} i \omega s^{\rho_j}} \; .
		\end{array} \right.
	\end{equation*}
\end{cor1}

\begin{proof}
	It suffices to note that $\phi_{n+1}^{(j)}(.,\omega,\rho_j,\mu_j)$ is the restriction to $(0,s_j] \subset D_j$ of the function $H_{n,\omega}^{(j)}\circ h$. Hence we remark that Proposition \ref{PROP5} affirms that $\phi_{1}^{(j)}(.,\omega,\rho_j,\mu_j) : (0,s_j] \longrightarrow \C$ is a primitive of $s \in (0,s_j] \longmapsto s^{\mu_j-1} e^{(-1)^{j+1} i \omega s^{\rho_j}}$, and use Proposition \ref{PROP6} to show that a primitive of $\phi_{n}^{(j)}(.,\omega,\rho_j,\mu_j) : (0,s_j] \longrightarrow \C$ is given by $\phi_{n+1}^{(j)}(.,\omega,\rho_j,\mu_j) : (0,s_j] \longrightarrow \C$, for $n \geqslant 1$.
\end{proof}

\begin{rem2} \label{REM2}
	\begin{enumerate}
		\em \item The function $\phi_{n+1}^{(j)}(.,\omega,\rho_j,\mu_j) : \, (0,s_j] \longrightarrow \C$ can be extended to $(0,+\infty)$. Indeed, we recall a parametrization of the curve $\Lambda^{(j)}(s)$ given by
		\begin{equation*}
			\lambda_{s}^{(j)} : t \in (0,+\infty) \longmapsto s+t e^{(-1)^{j+1} i\frac{\pi}{2 \rho_j}} \in \Lambda^{(j)}(s) \; ,
		\end{equation*}
		and we consider the following estimate
		\begin{equation} \label{eqrefe}
			\forall \,  t > 0 \qquad \Big|F_{n,\omega}^{(j)} \big( \lambda_{s}^{(j)}(t) , s \big) \Big| \leqslant \frac{1}{n!} \, t^{n+\mu_j-1} \, e^{-\omega t^{\rho_j}} \; ,
		\end{equation}
		which was obtained by noting that
		\begin{equation*}
			t \leqslant \left| s+t e^{(-1)^{j+1} i\frac{\pi}{2 \rho_j}} \right| = \left| \lambda_s^{(j)}(t) \right| \qquad \Longrightarrow \qquad t^{\mu_j-1} \geqslant \left| \lambda_s^{(j)}(t) \right|^{\mu_j-1} \; .
		\end{equation*}
		We notice that the right-hand side of \eqref{eqrefe} is an integrable function with respect to $t$ on $(0,+\infty)$, so $\phi_{n+1}^{(j)}(s,\omega,\rho_j,\mu_j)$ is well-defined for all $s > 0$.
		\item We can define the function $\phi_{n+1}^{(j)}(.,\omega,\rho_j,\mu_j) : \, (0,s_j] \longrightarrow \C$ at the point $0$. Indeed use estimate \eqref{eqrefe} one more time and remark that the right-hand side is independant on $s$. So by the dominated convergence Theorem, we can take the limit under the integral sign and we obtain
		\begin{equation*}
			\begin{aligned}
				\phi_{n+1}^{(j)}(0,\omega,\rho_j,\mu_j)	& := \lim_{s \rightarrow 0^+} \phi_{n+1}^{(j)}(s,\omega,\rho_j,\mu_j) \\
															& = \frac{(-1)^{n+1}}{n!} \int_{\Lambda^{(j)}(0)} z^{n+\mu_j -1} e^{(-1)^{j+1} i\omega z^{\rho_j}} dz \; .
			\end{aligned}
		\end{equation*}		
	\end{enumerate}
\end{rem2}


\begin{thebibliography}{99}
\bibitem{fam0} F. Ali Mehmeti, K. Ammari, S. Nicaise, \emph{Dispersive effects and high frequency behaviour for the Schrödinger equation in star-shaped networks, to appear in Port. Math.} \rm arXiv:1204.4998v2 [math.AP] (2014).
\bibitem{fam} F. Ali Mehmeti, R. Haller-Dintelmann, V. Régnier, \emph{The Influence of the Tunnel Effect on the $L^{\infty}$-time Decay}. Operator theory: Advances and Applications, \textbf{221} (2012), 11-24.
\bibitem{fam3} F. Ali Mehmeti, R. Haller-Dintelmann, V. Régnier, \emph{Energy Flow Above the Threshold of Tunnel Effect}. Operator Theory: Advances and Applications, \textbf{229} (2013), 65-76.
\bibitem{fam2} F. Ali Mehmeti, V. Régnier, \emph{Delayed reflection of the energy flow at a potential step for dispersive wave packets}. Math. Methods Appl. Sci., \textbf{27} (2004), 1145–1195.
\bibitem{banica} V. Banica, \emph{Dispersion and Strichartz inequalities for Schrödinger equations with singular coefficients}. SIAM J. Math. Anal, \textbf{35} (2003) no. 4, 868-838.
\bibitem{cazenave1998} T. Cazenave, F.B Weissler , \emph{Asymptotically self-similar global solutions of the nonlinear Schrödinger and heat equations}. Math. Z., \textbf{228} (1998), 83-120.
\bibitem{cazenave2010} T. Cazenave, J. Xie, L. Zhang, \emph{A note on decay rates for Schrödinger's equation}. Proceedings of the American Mathematical Society, \textbf{138} (2010) no. 1, 199-207.
\bibitem{erdelyi} A. Erdélyi, \emph{Asymptotics expansions}. Dover Publications, New York, 1956.
\bibitem{evans} L.C. Evans, \emph{Partial Differential Equations}. American Mathematical Society, USA, 1998.
\bibitem{hormander2} L. Hörmander, \emph{Propagation of singularities and semiglobal existence theorems for (pseudo)differential operators of principal type}, Ann. of Math. (2) \textbf{108} (1978), no 3, 569–609.
\bibitem{hormander} L. Hörmander, \emph{The Analysis of Linear Partial Differential Operators I}. Springer-Velag, Berlin Heidelberg New York Tokyo, 1983.
\bibitem{liess} O. Liess, \emph{Decay estimates for the solutions of the system of crystal optics}. Asymptotic Analysis, \textbf{4} (1991), 61-95.
\bibitem{msw} B. Marshall, W. Strauss, S. Wainger, \emph{$L^p-L^q$ estimates for the Klein-Gordon equation}. J. Math. pures et appl., \textbf{59} (1980), 417-440.
\bibitem{reed-simon} M. Reed, B. Simon, \emph{Methods of modern mathematical physics II : Fourier Analysis, Self-Adjointness}. Academics press, San Diego New York Boston London Sydney Tokyo Toronto, 1975.
\bibitem{strichartz} R.S. Strichartz, \emph{Restrictions of Fourier transforms to quadratic surfaces and decay of solutions of the wave equations}. Duke Math. J., \textbf{44} (1977), 705-714.
\end{thebibliography}
\end{document}